\documentclass[journal ]{new-aiaa}
\usepackage[utf8]{inputenc}
\usepackage{textcomp}

\usepackage{graphicx}
\usepackage{amsmath}
\usepackage[version=4]{mhchem}
\usepackage{siunitx}
\usepackage{longtable,tabularx}
\usepackage{caption}
\usepackage{subcaption}
\usepackage{algorithm2e}

\title{Rapid Determination of Low-Thrust Spacecraft Reachable Sets in Two-Body and Cislunar Problems}
\author{Sean R. Bowerfind \footnote{Ph.D. Student, Department of Aerospace Engineering, Auburn University, 141 Engineering Dr Auburn, AL 36849, AIAA Student Member.} and Ehsan Taheri\footnote{Assistant Professor, Department of Aerospace Engineering, Auburn University, 141 Engineering Dr Auburn, AL 36849, AIAA Senior Member.}}
\affil{Auburn University, Auburn, Al, 36830}

\begin{document}

\maketitle

\begin{abstract}
 The reachable set of controlled dynamical systems consist of the set of all possible reachable states from an initial condition, over a certain period of time under various control and operation constraints and exogenous disturbances. For space applications, determination of reachable sets is invaluable for trajectory planning, collision avoidance, ensuring safe and optimal performance in complex, real-world scenarios. Leveraging the connection between minimum-time and reachable sets, we propose a method for rapid determination of reachable sets for finite- and low-thrust spacecraft using an indirect minimum-time multi-stage formulation and the primer vector theory. 
 Reachable set analyses are presented for a minimum-time low-thrust Earth-Mars rendezvous problem and for several cislunar applications under circular restricted three-body dynamics. For the Earth-Mars problem and thruster parameters, results indicate that the minimum-time solution lies within the feasible set of position vectors, but on the boundary of the reachable set of velocity vectors. For the cislunar problems, L2 Halo orbit and Lunar Gateway 9:2 NRHO are considered. Our results indicate that the reachable set of low-thrust spacecraft coincides with invariant manifolds existing in the multi-body dynamical environments.
\end{abstract}
%Implementation details are presented for several problems. 
 %applied to a low-thrust spacecraft application . 
 %Using the thrust switching function of a minimum-fuel formulation, it is possible to estimate both feasible and reachable sets.  
 %a minimuWe use a minimum-fuel formulation to estimate both the reachable and feasible sets. 
 %Using a switching function, one can scale per-stage costate vectors to produce minimum-time trajectories (i.e., those that lie on the reachable set).
 %The ability to generate reachable sets is significant for analyzing system behavior, planning trajectories, verifying system properties, optimizing control strategies, and handling uncertainties. 
\section{Introduction}

Optimal control and robotics applications frequently use reachable set theory as a metric to assess cost and safety. In astrodynamics, solutions to reachability problems are useful for determining the spacecraft states that will lead to an inevitable collision or encounter, regardless of future collision avoidance maneuvers \cite{allen_machine_2014}. In the near Earth-orbit domain, the ability to accurately and rapidly compute the reachable set has significant implications in space-domain awareness to predict future spacecraft collisions and to generate optimal collision avoidance maneuvers. This is a problem currently facing all civil, commercial, and defense space-users. In 2022, the United States Space Force unveiled a new program under SpaceWERX titled ``Orbital Prime'', whose intentions are to stimulate the development of the Space Logistics/On-orbit servicing (OSAM/ISAM) industries. The primary focus of this initiative is to develop technologies to perform satellite re-location, de-orbiting, space debris removal, and remote and close inspection \cite{Spacewerx_2023}. Each of these tasks requires complex capabilities to map the Earth orbit environment. Current technology can detect more than 27,000 spacecraft and debris in orbit \cite{Garcia_2015}, however, a framework for modeling each of these objects and determining their reachable sets at future time horizons is crucial to accomplishing the OSAM/ISAM tasks. Knowledge of an object's reachable set can be used to perform evasive maneuvers (in the case of the International Space Station) or to perform orbital rendezvous and docking (to accomplish many of the OSAM/ISAM tasks) \cite{liu_analysis_2021}. Beyond low-Earth orbits, the ability to compute reachable sets is advantageous to cislunar \cite{vendl2021cislunar} and interplanetary mission design, as reachable sets provide mission designers with information regarding future planetary encounters from a given set of initial conditions \cite{chen_reachable_2019}. Future NASA deep-space missions with numerous planetary flybys could benefit substantially from knowledge about the reachable set of their respective spacecraft. For example, NASA's Europa Clipper, slotted to launch in October 2024, will orbit Jupiter and perform a projected 50 flybys of Europa \cite{NASA_2023}. With knowledge of a reachable set from an initial point in Jupiter's orbit, Europa Clipper can investigate low-thrust maneuvers that could encounter more of Jupiter's moons in addition to Europa to provide a higher scientific return on investment for the mission. 
 
%\subsection{Literature Review}

A reachable set is defined as the set of states that can be reached from a given initial condition  within a specified time span. Determining the reachable set is a problem that suffers from the ``curse of dimensionality'' since the Hamilton-Jacobi-Bellman (HJB) equation \cite{kirk_optimal_1970}, a partial differential equation (PDE), is traditionally used to determine reachable sets. The exact reachable set boundary can be found at the zero-level sets of the viscosity solution of an HJB problem \cite{crandall_properties_1984}. Across all classes of systems, linear, time-invariant (LTI) systems have been the most thoroughly analyzed. Results of LTI studies have shown that an analytical solution for the reachable set exists for cases of convex polyhedral state and control spaces \cite{borrelli_predictive_2017}. However, once constraints that lead to a nonconvex control space are introduced (common in astrodynamics problems), analytical solutions do not exist and numerical methods (i.e., \textit{direct} and \textit{indirect}) must be used instead \cite{trelat2012optimal}. Indirect methods will attempt to solve the viscosity solution of the HJB, while direct methods attempt to solve a (non-Hamiltonian) two-point boundary-value problem (TPBVP) \cite{wallace_applications_2005,chen_minimum-fuel_2021}. In 1980, Vinter showed how to approximate the reachability problem as a convex optimization problem and did the first study to demonstrate how to approximate the reachable set using a series of smoothing functions \cite{vinter_characterization_1980}. Determination of the reachable set becomes an increasingly difficult problem as the dimensionality of the problem and the number of control signals increases \cite{jiang_reachable_2023}. It rapidly becomes computationally infeasible to compute every possible trajectory variation through simulation alone \cite{mitchell_time-dependent_2005}. Modern research has expanded these approximation techniques to generate over- or under- approximations for continuous time-linear systems \cite{girard_efficient_2008}, general nonlinear dynamics \cite{goos_reachability_2003}, and more complex systems \cite{bansal2021deepreach}.

In the field of astrodynamics, reachability problems for both chemical and electric propulsion spacecraft are studied. Chemical propulsion reachability largely focuses on single impulse maneuvers and has been studied extensively for the cislunar environment \cite{xue_reachable_2010,zhang_reachable_2013,chen_new_2018,duan_simple_2019,wen_precise_2014,xuehua_reachable_2011,vinh_reachable_1995,wen_orbital_2014,lu_design_2021,wen_reachable_2023,wen_reachability_2019,komendera_intelligent_2012,chernick_closed-form_2021,xia_reachable_2022}.
Multiple electric propulsion or continuous thrust problems have been proposed for low-Earth, cislunar, and interplanetary mission planning applications \cite{lin_continuous-thrust_2023,chen_minimum-fuel_2021,pang_reachable_2022,lee_reachable_2018,wang_analytical_2023,jiang_reachable_2023}.
Other methods attempt to calculate reachable sets by utilizing passive $\Delta V$ maneuvers such as gravity assists and solar radiation pressure in small-body proximity operations. \cite{surovik_adaptive_2015,wen_calculating_2022,chen_reachable_2019, liu_analysis_2021,wu_reachable_2023}.

For the chemical propulsion case, an on-orbit range equation is derived for $\Delta V$ with J2 perturbations \cite{holzinger_-orbit_2014}. These results are used in \cite{chernick_closed-form_2021} to demonstrate closed-form impulsive control schemes by computing a minimum $\Delta V$ then solving a geometric path-planning problem. For the electric propulsion case, \cite{holzinger_reachability_2012} derives continuous indirect equations incorporating ellipsoid uncertainty with no mass loss. Several studies, \cite{kulumani_systematic_2019} and \cite{bando_nonlinear_2018} formulate optimal control problems and directly solve the HJB equation to generate reachable sets in the CR3BP over short time horizons. Specifically, \cite{kulumani_systematic_2019} takes a Poincar\'e section approach and details analysis with invariant manifolds and \cite{bando_nonlinear_2018} utilizes a power-series approximation of the HJB to determine periodic orbit transfer trajectories. A cross-disciplined  reachability-based trajectory design (RTD), developed by Kousik, has been used to perform safe path planning for quadcopters, wheeled-vehicles, and robotic manipulators \cite{kousik_bridging_2020}, \cite{kousik_technical_2019}. RTB computes the forward reachable set (FRS) of a vehicle offline using zonotope reachability \cite{bird_hybrid_2023}. The FRS is used in real time to perform obstacle avoidance maneuvering by comparing the FRS with an environment map that dictates the safe set. RTD uses this method to perform path planning in real time before feeding the safe trajectory to a low-level tracking controller. RTD has been successfully demonstrated in both simulation and hardware demonstrations.

%This paper focuses on reachable sets, however, the problem structure also provides a framework for finding the feasible set. Reachability implies feasibility in this case, since the feasible set consists of all points interior to the reachable set. 
%the convex or nonconvex multidimensional level sets of 
%This principle can be easily visualized when examining reachable position sets, in that any position state interior to the boundary of reachable positions can be considered a feasible position state. 
%The most important difference between reachable and feasible sets is the required control profile that is produced by the trajectories that make up each set. 
Determination of the reachable sets of low-thrust trajectories can be achieved by considering minimum-time solutions \cite{taheri_how_2020}. 
%Since the reachable set represents the ``limit case'' for each of the states in the considered problem, the control required to achieve a reachable state is also at its constraint. 
For a spacecraft equipped with low-thrust propulsion system, this means that the thruster will always be operating at maximum thrust and only the direction of the thrust vector is subject to variations across the  reachable set. 
%This observation is not consistent with trajectories in the feasible set. The only restriction placed on the trajectory in the problem formulation is that a minimum-fuel cost function is explored. As a result, this structure allows for bang-bang control profiles to achieve a minimum-fuel feasible trajectory. 
Reachable points are possible only by trajectories that use a thruster at maximum magnitude for the entire time horizon (i.e., a minimum-time trajectory). In fact, there is a duality between minimum-time and minimum-thrust solutions \cite{taheri_how_2020}. Trajectories that characterize feasible sets (i.e., interior to the reachable set) will produce a bang-bang control profile, correspond to minimum-fuel trajectories for a fixed time horizon and thrust parameters. 
%This algorithm uses an analytical switching function that can be used to obtain different bang-bang thrust profiles, and thus, different points on the reachable and feasible sets. However, 
The boundary of reachable set is defined by minimum-time solutions and  we only focus on the boundary of \textit{reachable} sets corresponding to minimum-time trajectories due to notable simplifications and speedup achieved combined witht the primer vector theory. 
%The scaling of the costates in made possible by the homogeneity property of the costate differential equations, switching function and the entire minimum-fuel low-thrust trajectory optimization \cite{bryson_applied_1975}.
%, so, we use the switching function to scale all the sampled terminal stage costates such that the propulsion system is always ON 
%(i.e., a minimum-time trajectory is achieved), and the final spacecraft state corresponds to a point on the reachable set only. 

Our main contributions are as follows: First, we propose a rapid indirect multi-stage formulation (IMF) for determining the reachable set of low-thrust spacecraft by leveraging the theoretical connection between minimum-time trajectories and reachable sets. This theoretical connection leads to noticeable algorithmic simplification and speedup compared to the minimum-fuel method proposed in \cite{patel_no_2023,patel_rapidly_2022}. We formulate all problems from a forward-reachable perspective that provides insights on the aforementioned spacecraft collision problems and celestial encounters.
More specifically, the discussion and presentation of the details of the IMF is given, which is an automated first-order approximation algorithm that is capable of sampling the reachable set. This algorithm requires no initial guess (automated) and utilizes only one first-order integration to forward propagate the nonlinear equations of motion for each sample. %The algorithm can be used to visualize either the reachable set, the boundary case for minimum-fuel trajectory optimization problems that is also a minimum-time solution, or the feasible set, which are interior points found inside the reachable boundary. 
%Additionally, the algorithm stores all sampled trajectories and can be used to present just the reachable set or a reachable tube that demonstrates the evolution of the reachable set through time. 
%Next, we show that the formulation of the reachable set algorithm allows for interchanging of different cost functions, as well as incorporation of initial boundary conditions on state uncertainty and impulse maneuvers. 
Secondly, we present the low-thrust reachability results for interplanetary trajectory optimization. In particular, minimum-time, rendezvous-type, Earth-to-Mars (Earth-Mars) maneuvers for a low-thrust spacecraft are presented. We show how the results can be used to rapidly determine position and velocity reachability subject to additional initial boundary constraints for a low-thrust Earth-Mars rendezvous trajectory optimization problem. In this context, the distinction between position- and velocity-reachability are highlighted. Thirdly, we applied the low-thrust reachability to a number of problems in the CR3BP of the Earth-Moon system. We demonstrate the versatility of this approach by generating novel results on the reachable set for a low-thrust spacecraft in cislunar space. We illustrate the unique evolution of the reachable set for a L2 Halo orbit and 9:2 NRHO. The last contribution of the paper is to offer theoretical insights into the behavior of the reachable sets in the CR3BP. We show that the reachable sets are connected to invariant manifolds of periodic orbits. The results also demonstrate the computational features of the proposed method.

The remainder of the paper is organized as follows. First, the IMF optimal control problem is introduced in Section \ref{sec:IMS}, which precedes discussions on the two-body equations of motion in Section \ref{sec:twobody}, the CR3BP in Section \ref{sec:cr3bp}, a derivation of minimum-time optimal control relations in Section \ref{sec:ocp}, and details on the application of initial boundary conditions in Section \ref{sec:boundarycons}. Next, we outline the rapid reachable set determination algorithm in Section \ref{sec:algorithm}. Finally, results are presented for the two-body dynamics in Section \ref{sec:twobodyresults} and for the CR3BP in Section \ref{sec:cr3bpresults}. We add remarks on invariant manifolds in Section \ref{sec:invariantmanifolds} to explain some trends identified within the CR3BP. Section \ref{sec:conclusion} presents a conclusion. In the remainder of the paper, vectors are denoted in boldface letters.

\section{Indirect Multi-Stage Formulation of Optimal Control Problems}
%We review the IMF, introduce the equations of motion, and derive the minimum-time optimal control equations that are foundational for algorithm implementation.

\subsection{Indirect Multi-Stage Formulation} \label{sec:IMS}
An IMF of OCPs is adopted, which offers notable  computational advantages since it allows for dynamics, constraints, and even cost functions to vary between different stages. This is an important feature for the introduction of initial boundary conditions, since incorporation of boundary conditions is achieved through a different cost function to account for any mass discrepancies due to an impulse maneuver. The other reason for using an IMF is the ability to derive analytical expressions for control, but also the fact that the entire reachable set determination is based on linearization.% compared to a direct method, which parameterizes all states and controls \cite{trelat_optimal_2012}. 

\begin{figure}[hbt!]
\centering
\includegraphics[scale = 0.2]{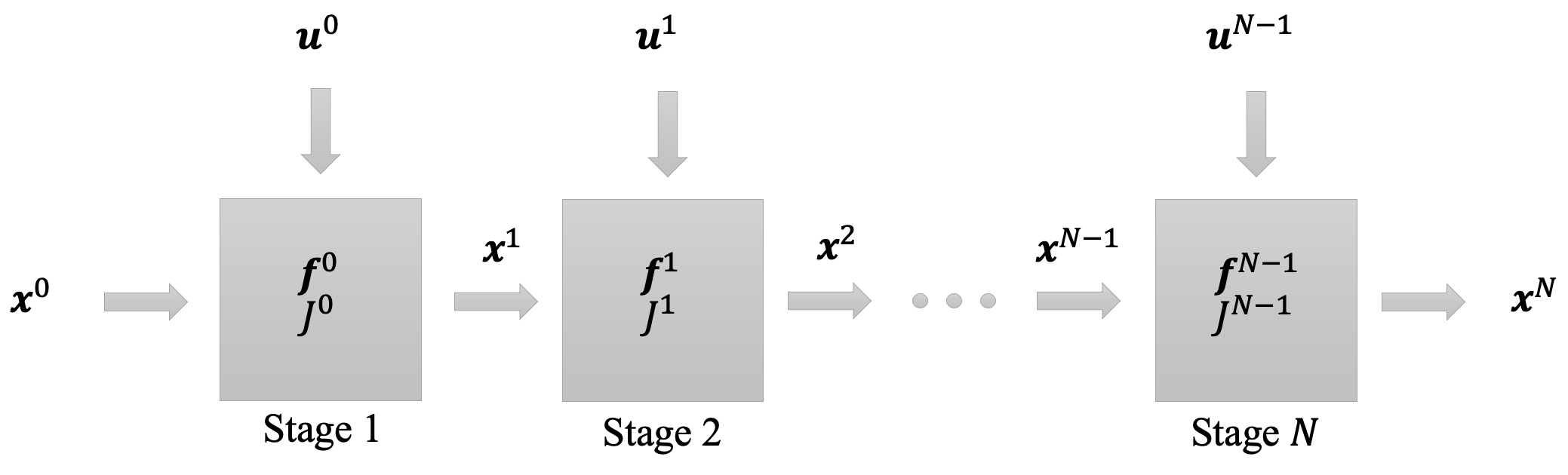}
\caption{Flowchart for a multi-stage formulation.}
\label{fig:multistageflowchart}
\end{figure}

The IMF divides a complete trajectory into a series of trajectory stages, but ``connects'' each stage through a sequence of equality constraints (see Fig. \ref{fig:multistageflowchart}). Each stage receives inputs of the states and controls from the previous stage, but independently optimizes control over only its respective stage \cite{bryson_applied_1975}. The entire trajectory can be discretized into $N$ stages and the $i^\text{th}$ stage (for $i\in\{1,\cdots,N\}$) can be modeled with its dynamical system consisting of states, $ \boldsymbol{x} \in \mathbb{R}^{n_x}$, controls, $ \boldsymbol{u} \in \mathbb{R}^{n_u}$, and parameters, $ \boldsymbol{p} \in \mathbb{R}^{n_p}$, where superscript is used to denote the stage index as,
\begin{align}
  \boldsymbol{x}^{i+1}& =\boldsymbol{F}^i\left(\boldsymbol{x}^i, t^i, \boldsymbol{u}^i; \boldsymbol{p}\right), & \text{with} &  & \boldsymbol{F}^i& =\boldsymbol{x}^i+\int_{t_i}^{t_{i+1}} \dot{\boldsymbol{x}} d t=\boldsymbol{x}^i+\int_{t_i}^{t_{i+1}} \boldsymbol{f}^{i-1}(\boldsymbol{x}, t, \boldsymbol{u}; \boldsymbol{p}) d t.
\end{align}

The IMF solves OCPs recursively, that is, stage $i=1$ is solved first, then $i=2,\cdots, i=N$ until the solution to all stages are computed. The derivation of the optimal control expressions is similar to continuous-time indirect formulations, except for some additional notation bookkeeping to track stages. The cost function for an IMF is:
\begin{equation}
    J=\phi\left(\boldsymbol{x}^N, N\right)+\sum_i L^i\left(t^i, \boldsymbol{x}^i, \boldsymbol{u}^i\right),
\end{equation}
where $\phi$ denotes a final cost and $L(\cdot)$ is the stage cost evaluated over the time span of each independent stage. The Hamiltonian of the complete trajectory can be written as the summation over the Hamiltonian of individual stages as,
\begin{align}
    H=\sum_i H^i,\quad \text{with} \quad  H^i= \left (\boldsymbol{\lambda}^{i+1} \right ) ^{\top}  \boldsymbol{F}^i+L^i+\sum_k v_k C_k.
    \label{hamcostate}
\end{align}

The costates, $\boldsymbol{\lambda} \in \mathbb{R}^{n_x}$, are specified with an $i+1$ index for convenience in the implementation of the algorithm, which includes a backward-in-time integration for the costates. In Eq.~\eqref{hamcostate}, $C_k$ and $v_k$ indicate the $k^\text{th}$ path constraint and its associated Lagrange multiplier (dimension of constraints, $k$ depends on the problem and the constraints). For each stage, partial derivatives of the Hamiltonian with respect to states and controls can be formed as,
\begin{align} \label{eq:costatemap}
    \frac{\partial H^i}{\partial \boldsymbol{x}} & = H_{\boldsymbol{x}}^i=\boldsymbol{\lambda}^i=(\boldsymbol{F_x}^{i})^{\top}\boldsymbol{\lambda}^{i+1}+L_{\boldsymbol{x}}^i+\sum_k v_k C_{k,\boldsymbol{x}}, & \frac{\partial H^i}{\partial \boldsymbol{u}} & = H_{\boldsymbol{u}}^i=\textbf{0}=(\boldsymbol{F_u}^{i})^{\top} \boldsymbol{\lambda}^{i+1}+L_{\boldsymbol{u}}^i+\sum_k v_k C_{k,\boldsymbol{u}}.
\end{align}

%According to the strong form of optimality, the partial derivative of the Hamiltonian with respect to control must be zero, which provides an extremal control law as,
%Note that Eq.~\eqref{hamcostate} has been set up to be integrated backward in time, which will be used later. 
%It will be shown how to derive the analytical optimal control expression for the problem of interest in the next section. Additionally, we will show how to incorporate initial boundary conditions using the IMF following the procedures developed in \cite{patel_no_2023,patel_rapidly_2022}.

\subsection{Two-Body Equations of Motion} \label{sec:twobody}

For the first results presented in this paper, the heliocentric phase of an  Earth-Mars rendezvous problem (with zero hyperbolic excess velocity) will be considered and third-body perturbations are ignored. Let $\boldsymbol{r} \in \mathbb{R}^3$ and $\boldsymbol{v}\in \mathbb{R}^3$ denote the position and velocity vectors, respectively. Let  $\boldsymbol{x} = [\boldsymbol{r}^{\top}, \boldsymbol{v}^{\top}]^{\top}$ denote the state vector and let $\boldsymbol{\hat{\alpha}}$ denote the thrust steering unit vector. the state dynamics can be written as \cite{bate_fundamentals_1971}.
\begin{equation}\dot{\boldsymbol{x}}=\left[\begin{array}{c}
\dot{\boldsymbol{r}} \\
\dot{\boldsymbol{v}}
\end{array}\right]=\left[\begin{array}{c}
\boldsymbol{v} \\
-\frac{\mu}{||\boldsymbol{r}||^3}\boldsymbol{r}+\frac{T_{\text{max}}}{m} \hat{\boldsymbol{\alpha}}\delta
\end{array}\right],
\end{equation}
where $\mu$ denotes the gravitational parameter of the Sun, $c = I_\text{sp} g_0$ is the net exhaust velocity, with $I_\text{sp}$ denoting the specific impulse, and $g_0$ being the gravitational acceleration at the Earth's surface. $T_{\text{max}}$ is the maximum thrust of the propulsion system operating in a thrust direction, $\hat{\boldsymbol{\alpha}}\in\mathbb{R}^3$, with a throttle setting of $\delta\in[0,1]$. In a minimum-fuel solution, the optimal control  of a low-thrust propulsion system will involve finding both $\hat{\boldsymbol{\alpha}}$ and $\delta$, however, for a minimum-time solution, it is known that the engine has to operate at maximum throttle during the entire trajectory ($\delta(t)=1~\forall t\in [t_0,t_f]$) \cite{taheri2017co}. Thus, different minimum-time trajectories can be generated by varying thrust directions, $\boldsymbol{\hat{\alpha}}$ all of which characterize the low-thrust reachable set. The spacecraft mass, $m$, can be calculated as, $ m(t) = m_0 - \frac{T_{\text{max}}}{c}(t-t_0)$.
%\begin{equation} \label{massflow}
%    m(t) = m_0 - \frac{T_{\text{max}}}{c}(t-t_0).
%\end{equation}
This minimum-time formulation allows the control vector, $\boldsymbol{u} \in \mathbb{R}^3$, to be stated simply as
\begin{align} \label{con}
    \boldsymbol{u}(t)&=
\hat{\boldsymbol{\alpha}}(t), & \Vert \hat{\boldsymbol{\alpha}}(t) \Vert& = 1
.
\end{align}

Alternative approaches for modeling a low-thrust spacecraft invoke a minimum-fuel formulation, which necessitates the tracking of mass and additional parameterization of the control vector by including both thrust magnitude and thrust direction ~\cite{junkins_exploration_2019}. A similar reachable set determination method proposed by Patel in Refs. \cite{patel_rapidly_2022,patel_no_2023,patel_reachable_2023, patel_rapid_2023} formulates a minimum-fuel problem and utilizes an inverse mass parameter, $n=1/m$, to track mass flow. This method  resolves some common problems encountered when using standard mass flow dynamics,  $\dot{m} = -\Vert\boldsymbol{T}\Vert/c$, in that this relation is non-differentiable when $\boldsymbol{T}=\boldsymbol{0}$. 

% The adopted parameterization of the control vector allows for the thrust vector, $\boldsymbol{T}$, to be linked to the position and velocity states, whereas the thrust magnitude, $T_\text{mag}$, is only used to track the change in mass (or inverse mass) of the spacecraft. Most popular implementations of Keplerian motion do not use inverse mass and simply track mass leading to the dynamics, $\dot{m} = -\Vert\boldsymbol{T}\Vert/c$. However, this formulation is non-differentiable when $\boldsymbol{T}=\boldsymbol{0}$, which is eliminated by using $n$ and separating $T_\text{mag}$. An alternative approach that is convenient for parameterizing the control vector is proposed in Ref. \cite{junkins_exploration_2019}. Two constraints will be imposed on the control that will simplify the formulation of the reachable set determination problem. Thrust magnitude must be less than the maximum capable by the propulsion system due to power limitations or thruster design considerations. 

% These constraints require that the maximum allowable thrust, $T_\text{max}$, and the thrust magnitude, $T_\text{mag}$, must be equal for a particular candidate solution to lie on the associated extremal solution. When imposed, these two constraints allow for a bang-bang control structure to be possible, which is typical of minimum-fuel low-thrust trajectories. In addition, when $T_\text{mag} = T_\text{max}$ along the entire trajectory, this corresponds to minimum-time maneuvers. Minimum-time solutions as will be shown correspond to the boundary of the feasible set (i.e., the  reachable set).

\subsection{Circular Restricted Three Body Problem} \label{sec:cr3bp}

The reachable set analysis of spacecraft in three-body dynamics has significant utility for near-term cislunar reachability applications. We consider the CR3BP model, which as will be shown offers more insights into the class of solutions that are obtained. That is, the problem utilizes the rotating synodic frame where the motion of the two primaries is assumed to be circular and the mass of the third object is negligible. The CR3BP is normalized such that the distance between the two primaries, total mass of the system, and mean motion of the two primary masses is unity. The characteristic length is defined as the distance between the two primaries, $l^*$, characteristic mass is defined as the sum of the masses of the two primaries, $m^* = m_1 + m_2$, and characteristic time is the mean motion of the primary system, $t^* = \sqrt{\frac{l^{*^3}}{Gm^*}}$, where $G$ is the universal gravitational constant. In the CR3BP, the mass ratio is an important parameter used to differentiate between different systems, defined as, $\mu = \frac{m_2}{m^*}$. These characteristic values for the Earth-Moon system are as follows \cite{thangavelu_transfers_2019}: $l^* = 3.844\times10^5$ km, $m^* = 6.0458\times10^{24}$ kg, $t^* = 375200$ seconds, and $\mu = 0.0121505856$.

%\begin{equation}
%    l^* = |\vec{r}_1 - %\vec{r}_2| .
%    \label{l*}
%\end{equation}

%The characteristic mass is defined as the sum of the masses of the two primaries

%\begin{equation}
%    m^* = m_1 + m_2.
%    \label{m*}
%\end{equation}

%The characteristic time is the mean motion of the primary system given as, 

%\begin{equation}
%    t^* = %\sqrt{\frac{l^{*^3}}{Gm^*}}, 
%    \label{t*}
%\end{equation}

%\begin{equation}
%    \mu = \frac{m_2}{m^*}
%    \label{mu}
%\end{equation}

%given in Table \ref{tab:charvals} 
%\begin{table}[]
%    \centering
%    \begin{tabular}{|c|c|c|}
%       \hline
%        Parameter & Value & Units \\
%        \hline
%        $l^*$ & $3.844\times10^5$ & km \\
%        \hline
%        $m^*$ & %$6.0458\times10^24$ & kg \\
%        \hline
%        $t^*$ & $375200 $& s \\
%        \hline
%        $\mu $& $0.0121505856$ & - \\
%        \hline
%    \end{tabular}
%    \caption{Characteristic Values for the Earth-Moon System \cite{thangavelu_transfers_2019}.}
%    \label{tab:charvals}
%\end{table}
In addition to the natural dynamics involving the motion of the two primaries, additional terms involving spacecraft mass and thrust must be normalized. Since the mass of the spacecraft was assumed negligible for the CR3BP derivation, these terms must be normalized by the initial mass of the spacecraft to avoid numerical issues. Since $m$ does not explicitly depend on the state, $\boldsymbol{x}$, we normalize control accordingly, $ m^*_{s/c}  = m_{0,s/c}$ and $T^*_{s/c} = \frac{m^*_{s/c} l^*}{{t^*}^2}$. Incorporating the geometry of the CR3BP presented in \cite{wiesel_spaceflight_1997}, the origin of the rotating frame is the barycenter of the two primaries and the $x$-axis is the vector collinear with the center of mass of the two primaries. Also, $m_1$ and $m_2$, are located at distances $-\mu$ and $1-\mu$ along the $x$-axis, respectively. Let $\hat{\boldsymbol{\alpha}} =  [\alpha_x,\alpha_y,\alpha_z]^{\top}$ denote the thrust steering unit vector. The equations of motion for a low-thrust spacecraft in normalized coordinates are
\begin{equation}\dot{\boldsymbol{x}}=\left[\begin{array}{c}
\dot{\boldsymbol{r}} \\
\dot{\boldsymbol{v}}
\end{array}\right]=\left[\begin{array}{c}
\boldsymbol{v} \\
x + 2v_y - \frac{(1-\mu)(x+\mu)}{r_1^3} + \frac{\mu(x+\mu-1)}{r_2^3} + \frac{T_{\text{max}}}{m}\alpha_x \\
y - 2v_x - \frac{(1-\mu)y}{r_1^3} - \frac{\mu y}{r_2^3} + \frac{T_{\text{max}}}{m} \alpha_y \\
 - \frac{(1-\mu)z}{r_1^3} - \frac{\mu z}{r_2^3} + \frac{T_{\text{max}}}{m} \alpha_z
\end{array}\right],
\label{cr3bp}
\end{equation}
where $r_1 = \sqrt{(x+\mu)^2 + y^2 + z^2}$ and $r_2 = \sqrt{(x+\mu-1)^2 + y^2 + z^2}$ denote the distances from the Earth and Moon, respectively. The state and control vectors and the control constraint $\left(\Vert\hat{\boldsymbol{\alpha}} \Vert = 1\right)$ are consistent for both the two-body and CR3BP dynamics.

\subsection{Formulation of Minimum-Time Optimal Control Problems} \label{sec:ocp}

%The goal of this paper is to produce minimum-time trajectories and map them onto the associated reachable set. To formulate the OCP associated with minimum-time trajectories, we follow the standard indirect formalism of optimal control theory. 
%The process is, identify a cost functional to be minimized, write the optimal control Hamiltonian associated with the considered cost functional and dynamics, form the partial derivative of the Hamiltonian with respect to states and controls to obtain the costate differential equations and the optimal control law, and enforce any constraints to complete the analytical control law solution \cite{kirk_optimal_1970}. 
Let $\Delta t$ denote the stage time interval, where $\Delta t$ is held constant across all stages. The objective function is,
\begin{equation}
    \text{minimize} \quad J=\sum_i \Delta t.
    \label{cost}
\end{equation}

The Hamiltonian for each stage from the single-stage dynamics, $\boldsymbol{F}^i$, objective function, $J$, and the constraint from Eq.~\eqref{con} can be written as,
\begin{equation}
    H^i=\Delta t+(\boldsymbol{F}^i)^{\top} \boldsymbol{\lambda}^{i+1}+\nu_1\left(\Vert\hat{\boldsymbol{\alpha}}^i\Vert-1\right),
    \label{hamiltonian}
\end{equation}
where $\nu_1$ is a Lagrange multiplier used for augmenting thrust constraints to the Hamiltonian. The value of $\nu_1$ must be determined as part of the solution procedure. Let $\boldsymbol{\lambda}^{\top} = [\boldsymbol{\lambda}_{\boldsymbol{x}}^{\top},\boldsymbol{\lambda}_{\boldsymbol{v}}^{\top}]$ denote the costate vector. The costate differential equations are derived by forming the partial derivative of the Hamiltonian in Eq.~\eqref{hamiltonian} with respect to the state vector, $\boldsymbol{x}$, written as (we use $\boldsymbol{F}_{\boldsymbol{x}}^{i,\top} = (\boldsymbol{F}_{\boldsymbol{x}}^{i})^{\top} $),
\begin{equation}
    \frac{\partial H^i}{\partial \boldsymbol{x}} = H_{\boldsymbol{x}}^i=\boldsymbol{\lambda}^i=\boldsymbol{F}_{\boldsymbol{x}}^{i,\top} \boldsymbol{\lambda}^{i+1}.
    \label{costate}
\end{equation}

Note that Eq.~\eqref{costate} is indeed set up to be integrated backward with respect to time, since the costate vector from the previous stage, $\boldsymbol \lambda^{i+1}$, is on the right-hand side. This backward-in-time step is related to the negative sign in $\dot{\boldsymbol{\lambda}} = - [\partial H/\partial \boldsymbol{x}]^{\top}$ for continuous-time relations obtained from applying the Euler-Lagrange equation in a standard indirect method. The control law is found from the partial derivative of Eq.~\eqref{hamiltonian} with respect to the control vector, $\boldsymbol{u}$, as,
\begin{equation}
    \frac{\partial H^i}{\partial \boldsymbol{u}} = H^i_{\boldsymbol{u}}=\textbf{0}=\boldsymbol{F}_{\boldsymbol{u}}^{\top} \boldsymbol{\lambda}^{i+1}+\nu_1\left(
\frac{\hat{\boldsymbol{\alpha}}^i}{\Vert\hat{\boldsymbol{\alpha}}^i\Vert}
\right).
\label{control}
\end{equation}

In Eqs.~\eqref{costate} and \eqref{control}, matrices $\boldsymbol{F}_{\boldsymbol{x}}$ and $\boldsymbol{F}_{\boldsymbol{u}}$ denote the state transition and sensitivity matrices, respectively, for linear systems \cite{junkins_analytical_2009}.  To evaluate $\boldsymbol{F}_{\boldsymbol{x}}$ and $\boldsymbol{F}_{\boldsymbol{u}}$ matrices, their matrix differential equations must be propagated through each stage to evaluate their values at the final time of each stage. Let $\Phi = \boldsymbol{F}_{\boldsymbol{x}}$ and $\Omega = \boldsymbol{F}_{\boldsymbol{u}}$, the matrix differential equations \cite{junkins_analytical_2009} can be written as,
\begin{align} \label{phidot}
    \dot\Phi(t,t_i) & = \mathbb{A}(t)\Phi(t,t_i), & \dot\Omega(t,t_i) & = \mathbb{A}(t)\Omega(t,t_i) + \mathbb{C}(t), %\label{omegadot}
\end{align}
where the matrices $\mathbb{A}(t)\in\mathbb{R}^{n_{\boldsymbol{x}}\times n_{\boldsymbol{x}}}$ and $\mathbb{C}(t)\in\mathbb{R}^{n_{\boldsymbol{x}}\times n_{\boldsymbol{u}}}$ can be derived as,
\begin{align}
    \mathbb{A}(t) & = \frac{\partial \dot{\boldsymbol{x}}}{\partial \boldsymbol{x}} =\left[\begin{array}{ccc}
\mathbf{0} & \boldsymbol{I} \\
\frac{\partial \dot{\boldsymbol{v}}}{\partial r} & \mathbf{0}
\end{array}\right],& 
     \mathbb{C}(t) & = \frac{\partial \dot{\boldsymbol{x}}}{\partial \boldsymbol u} = \left[\begin{array}{cccc}
\mathbf{0}\\
\operatorname{diag}(\frac{T_{\text{max}}}{m})
\end{array}\right].
\end{align}

The IMF notation can be updated and we can rewrite Eq.~\eqref{phidot} in an integral form as,
\begin{equation}
    \boldsymbol{F}_{\boldsymbol{x}}^i=\frac{\partial \boldsymbol{x}^{i+1}}{\partial \boldsymbol{x}^i}=\Phi\left(t^{i+1}, t^i\right)=\int_{t^i}^{t^{i+1}} \mathbb{A}(t) \Phi\left(t, t^i\right) d t,~ \boldsymbol{F}_{\boldsymbol{u}}^i=\frac{\partial \boldsymbol{x}^{i+1}}{\partial \boldsymbol{u}^i}=\Omega\left(t^{i+1}, t^i\right)=\int_{t^i}^{t^{i+1}} [\mathbb{A}(t) \Omega\left(t, t^i\right)+ \mathbb{C}(t)] d t,
\end{equation}
%\begin{equation}
%    \boldsymbol{F}_{\boldsymbol{u}}^i=\frac{\partial \boldsymbol{x}^{i+1}}{\partial \boldsymbol{u}^i}=\Omega\left(t^{i+1}, t^i\right)=\int_{t^i}^{t^{i+1}} [\mathbb{A}(t) \Omega\left(t, t^i\right)+ \mathbb{C}(t)] d t,
%\end{equation}
with the initial conditions defined as $\Phi(t^i,t^i) = \boldsymbol{I}$ and $\Omega(t^i,t^i) = \boldsymbol{0}$ (with $\boldsymbol{I}$ denoting the $n_x \times n_x$ identity matrix). The final step of the IMF is to solve Eq.~\eqref{control} for an optimal control law.
%, which incorporates a switching structure to account for the inequality constraint from Eq.~\eqref{con}. Define the thrust unit vector as $\hat{\boldsymbol{T}} =  \boldsymbol{T}/\Vert\boldsymbol{T}\Vert$, and $\boldsymbol{F_T}$ and $\boldsymbol{F}_{T_\text{mag}}$ as the rows of $\boldsymbol{F}_{\boldsymbol{u}}$ associated with $\boldsymbol{T}$ and $T_\text{mag}$. Extracting the top row of Eq.~\eqref{control}, which relates to the thrust norm constraint, Eq.~\eqref{con}, and 
Applying the strong form of optimality, such that Eq.~\eqref{control} is equal to zero, leads to
\begin{equation}
\boldsymbol{0}=\boldsymbol{F}_{\boldsymbol{u}}^{\top} \boldsymbol\lambda^{i+1}+\nu_1 \frac{\hat{\boldsymbol{\alpha}}^i}{\Vert\hat{\boldsymbol{\alpha}}^i\Vert}, \rightarrow -\frac{\boldsymbol{F}_{\boldsymbol{u}}^{\top} \boldsymbol\lambda^{i+1}}{\nu_1}=\frac{\hat{\boldsymbol{\alpha}}^i}{\Vert\hat{\boldsymbol{\alpha}}^i\Vert}=\hat{\boldsymbol{\alpha}}^i.
\end{equation}

Using the definition of $\hat{\boldsymbol{\alpha}}$, $\nu_1$ must be equal to $||\boldsymbol{F}_{\boldsymbol{u}}^{\top} \boldsymbol\lambda^{i+1}||$, which yields the familiar primer vector optimal control law,
\begin{align} \label{law}
    \hat{\boldsymbol{\alpha}}^i& =-\frac{\boldsymbol{F}_{\boldsymbol{u}}^{i, \top} \boldsymbol{\lambda}^{i+1}}{\left\Vert \boldsymbol{F}_{\boldsymbol{u}}^{i, \top} \boldsymbol{\lambda}^{i+1}\right\Vert}.
\end{align}

\subsection{Initial Boundary Conditions} \label{sec:boundarycons}

The IMF allows for imposing constraints on only one stage, or even using a different cost functional across different stages. We take advantage of this feature to present the framework for imposing position and velocity ellipsoid uncertainty constraints (per the method of \cite{holzinger_reachability_2012}) as well as an impulse maneuver at the simulation start time, adapting the work presented in \cite{patel_rapidly_2022,patel_rapid_2023}. However, we present related numerical results to the Earth-Mars problem. To accomplish this, we add a pseudo-zero stage that converts the reference trajectory initial states,  $\boldsymbol{r}^\text{ref}$ and $\boldsymbol{v}^\text{ref}$, into perturbed initial states,  $\boldsymbol{r}^\text{*}$ and $\boldsymbol{v}^\text{*}$. We use a matrix-mapping technique to account for the variation in each of the states as,
\begin{equation} \label{eq:uncere}
    \left[\begin{array}{c}
\boldsymbol{r}^{*} \\
\boldsymbol{v}^*
\end{array}\right]=\left[\begin{array}{c}
\boldsymbol{r}^\text{ref} \\
\boldsymbol{v}^\text{ref}
\end{array}\right]+\left[\begin{array}{llll}
\boldsymbol{I} & \mathbf{0} & \mathbf{0}  \\
\mathbf{0} & \boldsymbol{I} & \boldsymbol{I} 
\end{array}\right]\left[\begin{array}{c}
\boldsymbol{\delta {r}} \\
\boldsymbol{\delta v_1}\\
\boldsymbol{\delta v_2}
\end{array}\right].
\end{equation}

Our goal is to solve for $\left[\boldsymbol{\delta r}, \boldsymbol{\delta v_1}, \boldsymbol{\delta v_2}\right]^\top$, which are the states associated with variations from the reference state. In Eq.~\eqref{eq:uncere}, $\boldsymbol{\delta r}$ and $\boldsymbol{\delta v_1}$ are associated with ellipsoid position and velocity uncertainty constraints. In addition, $\Delta V$ impulse maneuver constraints are handled through $\boldsymbol{\delta v_2}$ . An additional cost function is required to account for the variation in initial mass due to the impulse maneuver. We select Tsiolkovsky rocket equation as the cost function to track the mass consumed due to this impulsive maneuver as, $L^{*}=m_0-m_0 e^{-\frac{\Delta V}{c}}$. The (pseudo-zero stage) Hamiltonian  is formed as,
\begin{equation}
    \begin{aligned}
& H^{*}=m_0-m_0 e^{-\frac{\Delta V}{c}}+\boldsymbol\lambda^\top [\boldsymbol{r}^{*\top},\boldsymbol{v}^{*\top}]^{\top} +\sum_k \nu_k C_k. 
\end{aligned}
\end{equation}
%\left(\left[\begin{array}{c}
%r^\text{ref} \\
%v^\text{ref} \\
%n^\text{ref}
%\end{array}\right]+\left[\begin{array}{llll}
%\boldsymbol{I} & \mathbf{0} & \mathbf{0} & \mathbf{0} \\
%\mathbf{0} & \boldsymbol{I} & \boldsymbol{I} & \mathbf{0} \\
%0 & 0 & 0 & 1
%\end{array}\right]\left[\begin{array}{c}
%\boldsymbol{\delta  r} \\
%\boldsymbol{\delta  v_1} \\
%\boldsymbol{\delta  v_2} \\
%\delta n
%\end{array}\right]\right)

Next, we must express $C_k$ for both types of constraints. We use ellipsoid uncertainty constraint formulations from \cite{holzinger_reachability_2012} to express the position and velocity uncertainty in the form of equality constraints as,
\begin{align}
    \frac{1}{2} \boldsymbol{\delta} \boldsymbol{r}^{\top} \mathbb{E}_{\boldsymbol{r}} \boldsymbol{\delta} \boldsymbol{r}-\frac{1}{2} r_\text{ref}^2& =0,& \frac{1}{2} \boldsymbol{\delta} \boldsymbol{v}_{\mathbf{1}}^{\top} \mathbb{E}_{\boldsymbol{v}} \boldsymbol{\delta} \boldsymbol{v}_{\mathbf{1}}-\frac{1}{2} v_\text{ref}^2& =0,
    \label{posellip}
\end{align}
where $\mathbb{E}_{\boldsymbol{r}}$ and $\mathbb{E}_{\boldsymbol{v}}$ are matrices determining the major axes of the uncertainty ellipsoid.
We need expressions for $\boldsymbol{\delta r}$ and $\boldsymbol{\delta v_1}$ to successfully implement  constraints in Eq.~\eqref{posellip}. Computing the partial derivative of $H^*$ with respect to $\boldsymbol{\delta r}$ gives
\begin{equation} \label{dr}
    \frac{\partial H^*}{\partial\boldsymbol{\delta r}} = H_{\boldsymbol{\delta r}}^*=\textbf{0}=\boldsymbol\lambda_{\boldsymbol{r}}+\nu_1 \mathbb{E}_{\boldsymbol{r}} \boldsymbol{\delta r},  \rightarrow \boldsymbol{\delta r} =-(\nu_1 \mathbb{E}_{\boldsymbol{r}})^{-1} \boldsymbol\lambda_{\boldsymbol{r}}, 
\end{equation}
where $\boldsymbol\lambda_{\boldsymbol{r}}$ is the pseudo-zero stage  costates associated with the position vector. Substituting Eq.~\eqref{dr} into the position ellipsoid constraint, Eq.~\eqref{posellip}, yields an expression for $\nu_1$, which can be written as,
\begin{align}
    0  =\frac{1}{2} r_\text{ref}^2-\frac{1}{2} \frac{\boldsymbol\lambda_r^\top \mathbb{E}_{\boldsymbol{r}}^{-1} \boldsymbol\lambda_r}{\nu_1^2}, \rightarrow \nu_1 =\frac{\boldsymbol\lambda_r^\top \mathbb{E}_{\boldsymbol{r}}^{-1} \boldsymbol\lambda_r}{r_\text{ref}^2}.
\end{align}

Substituting $v_1$ into Eq.~\eqref{dr} gives the final equation as,
\begin{equation}
    \boldsymbol{\delta r} =-r_\text{ref} \frac{\mathbb{E}_{\boldsymbol{r}}^{-1} \boldsymbol\lambda_{\boldsymbol{r}}}{\sqrt{\boldsymbol\lambda_{\boldsymbol{r}}^\top \mathbb{E}_{\boldsymbol{r}}^{-1} \boldsymbol\lambda_{\boldsymbol{r}}}}.
    \label{poscon}
\end{equation}

All parameters in Eq.~\eqref{poscon} are given, aided by the backwards integration in the second stage of the algorithm, which yields the costates, $\boldsymbol{\lambda_r}$. The initial condition position uncertainty constraint can be implemented. It is possible to follow the same process to obtain a solution for $\boldsymbol{\delta v_1}$ as,
\begin{align}
    H_{\boldsymbol{\delta v_1}}^* =\textbf{0}=&  \boldsymbol{\lambda_v} +\nu_2 \mathbb{E}_{\boldsymbol{v}} \boldsymbol{\delta v_1},& \boldsymbol{\delta v_1}  =& -v_\text{ref} \frac{\mathbb{E}_{\boldsymbol{v}}^{-1} \boldsymbol{\lambda_v}}{\sqrt{\boldsymbol{\lambda_v}^\top \mathbb{E}_{\boldsymbol{v}}^{-1} \boldsymbol{\lambda_v}}}.
    \label{velcon}
\end{align}

Equations \eqref{poscon} and \eqref{velcon} are analytical closed-form expressions for the change in initial conditions due to an ellipsoid uncertainty. These changes can be implemented in the reachable set algorithm without any structural modification to the algorithm (see Eq.~\eqref{eq:uncere}). The other type of scenario is the $\Delta V$ impulse maneuver, which involves a separate boundary condition, since one must allow the magnitude of the $\Delta V$ impulse to take values less than or equal to the maximum capability of the propulsion system, $\Delta V_{\max }$. This constraint associated with the initial impulse can be written as, \begin{align} \label{impulsecons}
    \frac{1}{2} \boldsymbol{\delta v_2}^\top \boldsymbol{\delta v_2}-\frac{1}{2} \Delta V_{\max }^2 & \leq 0.
\end{align}

The derivation of the closed-form solutions for these constraints is more involved due to the introduction of an inequality constraint that can be active or inactive. A complete derivation of the solution for  $\boldsymbol{\delta v_2}$ and $\Delta V$ is given in the Appendix. However, if the inequality constraint is active, then we have the following two relations
\begin{align}\label{dv2bind}
    \Delta V & = \Delta V_{\max}, & \boldsymbol{\delta v_2} & = -\Delta V_{\max} \frac{\boldsymbol{\lambda_v}}{\left\Vert\boldsymbol{\lambda_v}\right\Vert}.
\end{align}

If the inequality constraint is not active, the optimal solution can be derived as,
\begin{align}\label{dv2nobind}
    \Delta V & =-c \ln \left(\frac{\left\Vert\boldsymbol{\lambda}_{\boldsymbol{v}}\right\Vert c+\sqrt{\left\Vert\boldsymbol{\lambda}_{\boldsymbol{v}}\right\Vert^2 c^2}}{2 m_0}\right), &  \boldsymbol{\delta} \boldsymbol{v}_{\boldsymbol{2}}& =c \ln \left(\frac{\left\Vert\boldsymbol{\lambda}_{\boldsymbol{v}}\right\Vert c+\sqrt{\left\Vert\boldsymbol{\lambda}_{\boldsymbol{v}}\right\Vert^2 c^2}}{2 m_0}\right) \frac{\lambda_{\boldsymbol{v}}}{\left\Vert\boldsymbol{\lambda}_{\boldsymbol{v}}\right\Vert}.
\end{align}

The final step is to address the initial mass perturbation due to an impulse maneuver, which can be written as, $\delta m = m_0 - m_0e^{-\frac{\Delta V}{c}}$
%\begin{equation} \label{dm}
%    \delta m = m_0 - m_0e^{-\frac{\Delta V}{c}},
%\end{equation}
with $\Delta V$ becoming the computed value using the logic presented above. Finally, we need to update the initial mass as, $ m^* = m_0 + \delta m$. We emphasize that the analytical solutions for an initial condition variation due to an impulse maneuver are only functions of the costates and other constants specified by the problem. This makes implementation of these types of constraints simple given the algorithm structure (i.e, treating stages differently).

\section{Reachable Set Determination Algorithm} \label{sec:algorithm}
This section presents details and schematics for the implementation of the reachable set algorithm, which consists of three main phases. First, a zero-thrust reference trajectory is used to obtain all reference states, state transition matrices, and sensitivity matrices corresponding to the reference trajectory. Second, a sufficiently large number of costates (at the final stage, $N$)  are sampled from a 6-dimensional unit ball corresponding to the costates associated with the position and velocity vectors. For each sampled costate vector, the costates are integrated backward in time, using the costate update mapping in Eq.~\eqref{costate}, until the initial stage. Using the costates, control laws given by Eq.~\eqref{law} are computed and stored. Last, the reachable set is reconstructed using a forward-in-time integration of the fully \textit{nonlinear} state dynamics with the computed optimal control. 
%These three steps are repeated for every sampled trajectory to produce a number of sample points forming the reachable set. 
The sampling-based aspect of the algorithm allows for thousands of reachable trajectories to be generated rapidly and in a parallel manner. Note that for each trajectory  there is only one numerical integration and no optimization solver is used (as opposed to a direct optimal control method). 

\begin{figure}[hbt!]
\centering
\includegraphics[scale = 0.43]{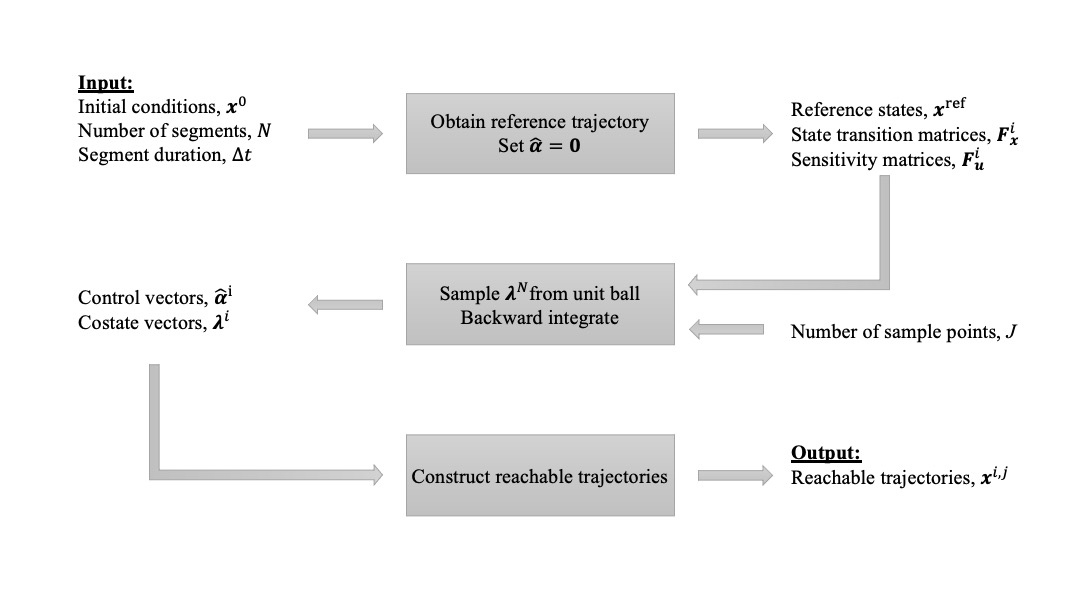}
\vspace{-10mm}
\caption{Flowchart for the reachable set determination algorithm.}
\label{flowchart}
\end{figure}

When selecting a reference trajectory, we select a zero-thrust (i.e., ballistic) trajectory. This is due to a trivial solution of the OCP in which $\boldsymbol{\lambda}^N=\textbf{0}$, causing $\{\boldsymbol{u}_i\}_i^N=\textbf{0}$. 
%To show this, insert $\boldsymbol{\lambda}^N=\textbf{0}$ into the switching function, Eq.~\eqref{sf}, which will return a negative value. The switching logic means that since $SF<0$, the thrust values will also be zero. To verify this is true for all stages,  propagate the costate differential equations backward in time using Eq.~\eqref{costate} to observe that all previous $i^\text{th}$ stage costates are zero, leading to an unpowered reference trajectory. 
Integrate the system dynamics forward-in-time from initial conditions $\boldsymbol{x}^{0}$ with no control to obtain the reference trajectory, $\boldsymbol{x}^{\text{ref}}$, and $\{\boldsymbol{F_x}^i\}_i^N$ and $\{\boldsymbol{F_u}^i\}_i^N$ matrices.

Next, we sample terminal stage costates $\boldsymbol{\lambda}^N$ from a 6-dimensional unit ball (i.e., $\left\Vert\boldsymbol{\lambda}^N\right\Vert=1$). The reason we uniformly sample $\boldsymbol{\lambda}^N$ is due to the control law in Eq.~\eqref{law}, which indicates that the magnitude of $\boldsymbol{\lambda}^N$ does not matter (the control law is homogenous in costates); hence, a reasonable method to determine these values is to uniformly sample them from a unit ball. All previous stage costates are recovered, backwards in time, from a recursive matrix multiplication outlined in Eq. \ref{costate}.
% It must be noted that simply sweeping over a unit ball only provides near-optimal control laws. The true optimality must be recovered by scaling the thrust magnitude correctly.
%The final stage is to recover all reachable trajectories.
%by correctly scaling the costates to ensure the switching function is positive. This is completed by choosing the maximum single-stage scaling factor for each sampled trajectory and multiplying all stage costate vectors by the maximum scaling factor. 
The optimal control can be computed from Eq.~\eqref{law} and a final forward-in-time integration is used to produce the entire reachable set.  The indices $i$ and $j$ in Fig.~\ref{flowchart} indicate the $i^{\text{th}}$-stage variable quantity and the $j^{\text{th}}$ trajectory forming the reachable set. The proposed minimum-time reachable set determination algorithm is outlined in Algorithm \ref{alg:two}.

\RestyleAlgo{ruled}

%% This is needed if you want to add comments in
%% your algorithm with \Comment
\SetKwComment{Comment}{/* }{ */}

\begin{algorithm}[hbt!]
\caption{Reachable Set Determination}\label{alg:two}
% \KwData{intial conditions, segment time interval, number of segments, number of sample points}
\KwResult{$\boldsymbol{x}^{i,j}$ Reachable Trajectories}
$\boldsymbol{x}^{0} \gets$ Initial Conditions\;
$\boldsymbol{p} \gets$ Constant Parameters\;
$N \gets$ Number of Segments\;
$\Delta t \gets$ Segment Duration\;
$J \gets $Number of sampled reachable points\;
$\boldsymbol{u}^{\text{ref}} \gets \boldsymbol{0}$\;
$t^i \gets 0$\;
\While{$i=0; i < N; i \gets i+1$}{
    $\boldsymbol{x}^{i+1},  \boldsymbol{F}_{\boldsymbol{x}}^i,  \boldsymbol{F}_{\boldsymbol{u}}^i \gets \boldsymbol{F}^i\left(\boldsymbol{x}^i, \boldsymbol{u}^{\text{ref}}, t^i, t^i+\Delta t, \boldsymbol{p}\right)$\Comment*[r]{Unpowered reference trajectory}
    $\boldsymbol{x}^{i}\gets\boldsymbol{x}^{i+1}$\;
    $t^i\gets t^i + \Delta t$\;
}
\For{$j=1:J$}{
    $\boldsymbol{\lambda}^N\gets$ Sampled from Unit Ball\;
    $i\gets N$\;
    \While{$i=N; i \geq 1; i\gets i-1$}{
    $\hat{\boldsymbol{\alpha}}^i\gets-\frac{\boldsymbol{F}_{\boldsymbol{u}}^{i, \top} \boldsymbol{\lambda}^{i+1}}{\left\Vert \boldsymbol{F}_{\boldsymbol{u}}^{i, \top} \boldsymbol{\lambda}^{i+1}\right\Vert}$\Comment*[r]{Optimal unit thrust steering vector}
    $\boldsymbol{\lambda}^{i}\gets\boldsymbol{F}_{\boldsymbol{x}}^{i^{\top}} \boldsymbol{\lambda}^{i+1}$\;
    }
}
\For{$j=1:J$}{
    \While{$i=0; i < N, i\gets i+1$}{
        \If{$i=0$}{
                $\boldsymbol{\delta r} \gets$ Eq.~\eqref{dr}\;
               $\boldsymbol{\delta v}_1 \gets$ Eq.~\eqref{velcon}\;
                $\Delta V \gets \Delta V_{\max}$\;
                \eIf{\text{Eq.}~\eqref{nu3_1}$>0$ }{
                        $\boldsymbol{\delta v}_2, \Delta V \gets$ Eq.~\eqref{dv2bind}\;
                }{
                        $\boldsymbol{\delta v}_2, \Delta V \gets$ Eq.~\eqref{dv2nobind}\;
                }
              $\boldsymbol{x}^{i} \gets$ Eq.~\eqref{eq:uncere}  \;
        }
    $\boldsymbol{x}^{i+1,j} \gets \boldsymbol{F}^i\left(\boldsymbol{x}^i, \hat{\boldsymbol{\alpha}}^{i,j}, t^i, t^i+\Delta t, \boldsymbol{p} \right)$\Comment*[r]{Reconstruct reachable trajectories}
    $\boldsymbol{x}^{i}\gets\boldsymbol{x}^{i+1}$\;
    $t^i\gets t^i + \Delta t$\;
    }
    $j\gets j+1$\;
}

\end{algorithm}

The proposed algorithm shares multiple aspects with a recent minimum-time based algorithm by Patel in \cite{patel_rapid_2023}, but we briefly discuss some key differences between the method dubbed ``FFOE'' and our proposed algorithm. First, FFOE utilizes a cost function designed to maximize final state in a unit-ball sampled direction, whereas, our proposed method formulates a minimum-time IMF problem that allows direct comparison with the solution to a minimum-time OCP (see Sec. \ref{sec:resultscompareOCP}). Additionally, our formulation simplifies the dynamics by removing mass and its costate differential equations from the set of state-costate vector. This is aided by the minimum-time formation, since it is known that mass can be integrated  analytically using thrust and time of flight \cite{taheri2017co}. Ref. \cite{patel_rapid_2023} tracks inverse mass, which increases the dimensionality of the FFOE and leads to a longer computational time compared to the proposed algorithm (independent of the platform and programming language). Last, \cite{patel_rapid_2023} formulates an impulse maneuver initial constraint by assuming that the applied $\Delta V$ is equivalent to $\Delta V_{\max}$ of the spacecraft. Our proposed method allows $\Delta V$ taking any value less than or equal to $\Delta V_{\max}$ at the expense of additional logic in the algorithm (see Appendix \ref{sec:appendix})

Nevertheless, we can use the combined observations from this study and \cite{patel_rapid_2023} to draw conclusions on the shortcomings of a ``sampling-type'' reachable set estimation method.  In general, the algorithm accuracy begins to deteriorate over long time horizons \textit{unless} the number of sample points increases. Additionally, since the core of this method is a linearization about an unpowered reference trajectory, we expect that as the ratio of the acceleration produced by the propulsion to the natural acceleration due to gravity of the gravitational body(s) increases, the accuracy of the reachable set algorithm will decrease. These issues were also identified in \cite{patel_rapid_2023} in which FFOE was used in an iterative manner where each sample trajectory was linearized about and the reachable set recomputed, however this method exhibited diminishing returns and numerical stability issues.

\section{Two-Body Dynamics Results} \label{sec:twobodyresults}

The results that have been obtained for the rapid reachability and rendezvous determination problem for a minimum-time Earth-Mars transfer are presented in this section. This example is taken from \cite{taheri_enhanced_2016}. For all results, $T_\text{max} =  0.5$ Newton, $I_\text{sp} =3000$ seconds, and $m_0 = 1000$ kg. The initial position and velocity vectors are $\boldsymbol{r} = [-140699693;-51614428; 980]^{\top}$ km and $\boldsymbol{v} = [9.774596;-28.07828;4.337725\text{E-4}]^{\top}$ km/s, respectively, 
%state vector is given as, $\boldsymbol x = [-140699693;-51614428; 980; 9.774596;-28.07828;4.337725\text{E-4}]^{\top}$, with units km, km/s
corresponding to a departure Epoch of April $10^{\text{th}}$, 2007. Segment duration, $\Delta t$ is set to 86400 seconds. The reachable set algorithm was implemented in MATLAB on a 2018 MacBook Pro with a quad-core 2.3 GHz processor using parallel computing capabilities with 4 CPU cores. Additionally, MATLAB's integrated code compiler was used to generate MEX files to speed up the algorithm. NASA's spice toolkit was used for accessing DE440 planetary ephemerides. Computations are completed in the J2000 ecliptic reference frame centered at the solar system barycenter. %In all simulation, $\Delta t = 1$ day is considered.

\subsection{Reachable set analysis for a fixed time horizon}\label{sec:SRSA}
\begin{figure}[!tbp]
\centering 
\subfloat[2D view]{%
  \includegraphics[clip,width=0.95\columnwidth]{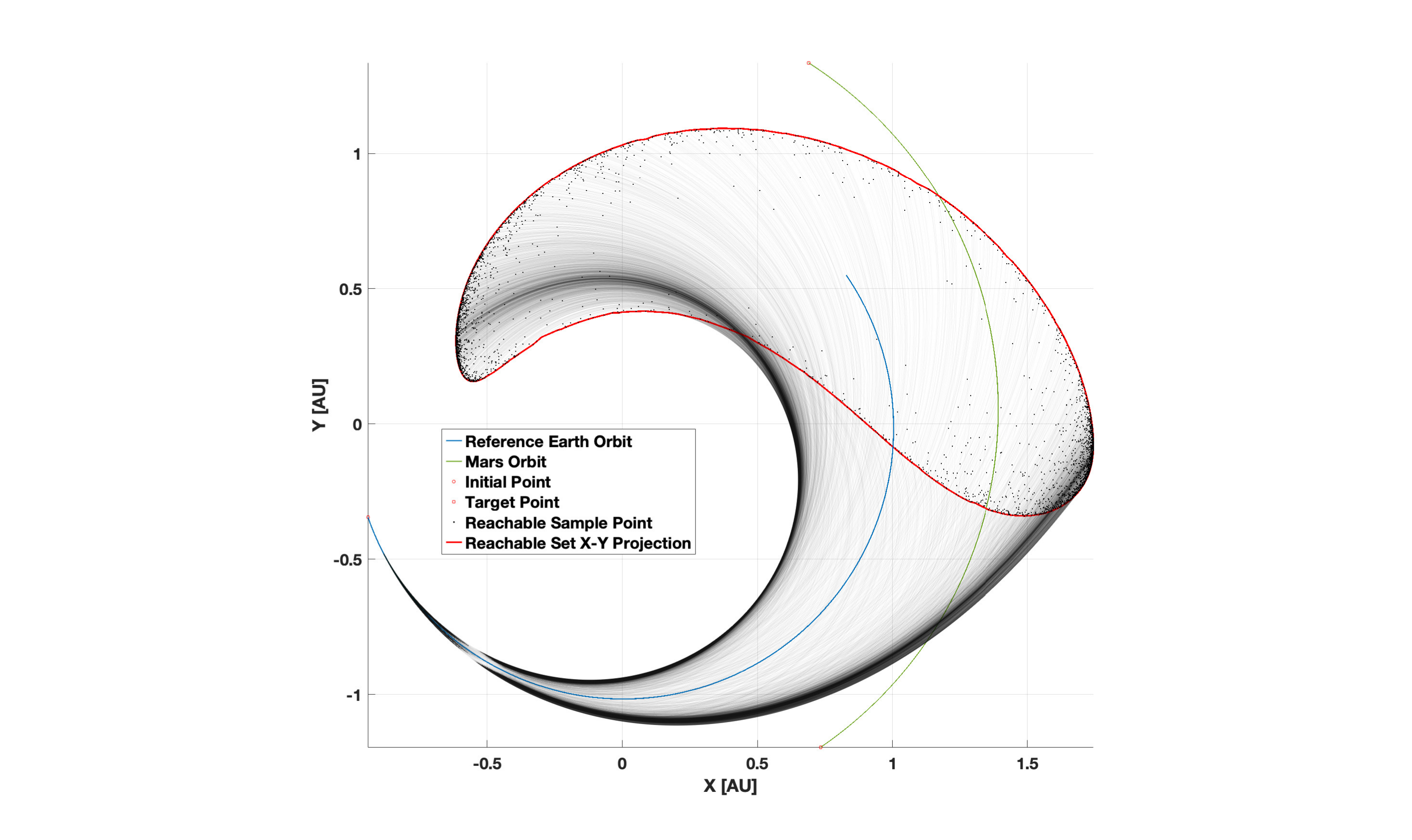}
  \label{fig:reacha}
}

\subfloat[3D view]{%
  \includegraphics[clip,width=0.95\columnwidth]{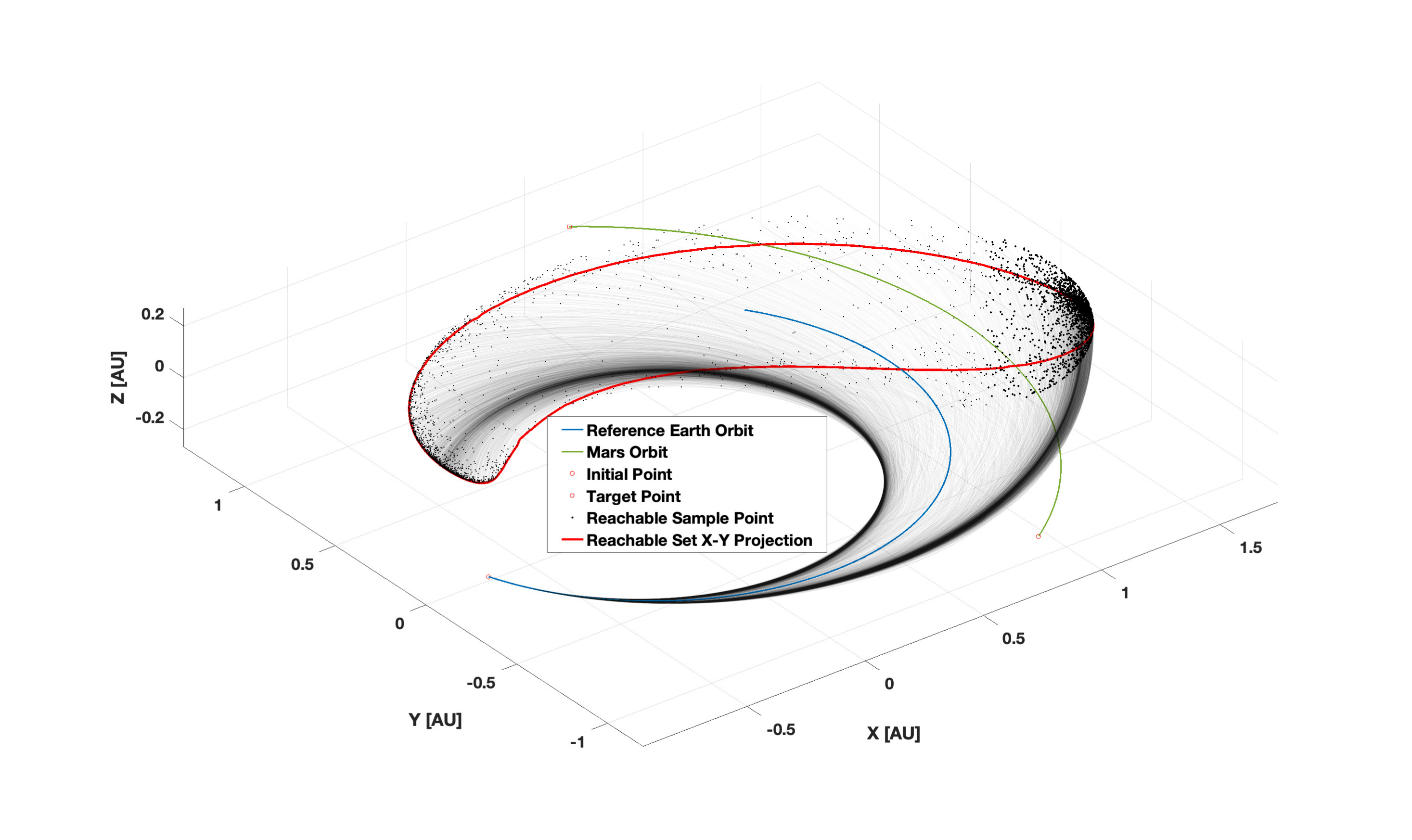}
  \label{fig:reachb}
}

\caption{Position reachable set over a 200-day time horizon with 5000 sample trajectories.}
\label{fig:EMreachabilitymain}
\end{figure}

%The reachable set algorithm was used to compute the reachable set for a low-thrust spacecraft with the given initial conditions for a 
We consider a time horizon of 200 days with a total of 5000 sample trajectories. This simulation completed in 13.61 seconds, which demonstrates rapid trajectory generation ability with 2.72 milliseconds per trajectory. The first, and most important, observation from Fig. \ref{fig:reacha} is that for the time span of 200 days, a portion of the Mars orbit is within the reachable set. However, the position of Mars (red square) is not within the reachable set (red). A longer time of flight is required for Mars to become reachable. The reference trajectory (blue line) is equivalent to the orbit of the Earth, since no control is applied. In Fig. \ref{fig:reacha}, AU stands for the astronomical unit (AU = $149 \times 10^6$ km).

Figure \ref{fig:reachb} depicts a three-dimensional (3D) view of the reachable set, which is a deformed ellipsoid-like 3D shape created by the endpoint of all reachable trajectories. A funnel-like structure exists as we plot the time history of all reachable trajectories (light gray lines), we refer to this as the reachable funnel. The boundary points of the reachable funnel act as a time history of all reachable sets within a 200-day time span. The feasible set and feasible funnel are all the position states classified as the interior of the reachable set and reachable funnel. As the time of flight increases, the number of sample trajectories has to be increased because for long-time-horizon maneuvers the volume of the reachable set expands rapidly, so more sample points are required to construct an accurate approximation of the reachable set.

\subsection{Reachable set comparison with the solution from a direct formulation}\label{sec:resultscompareOCP}

\begin{figure}[!tbp]
\centering 
\subfloat[Minimum-time trajectory and thrust vector.]{%
  \includegraphics[clip,width=0.9\columnwidth]{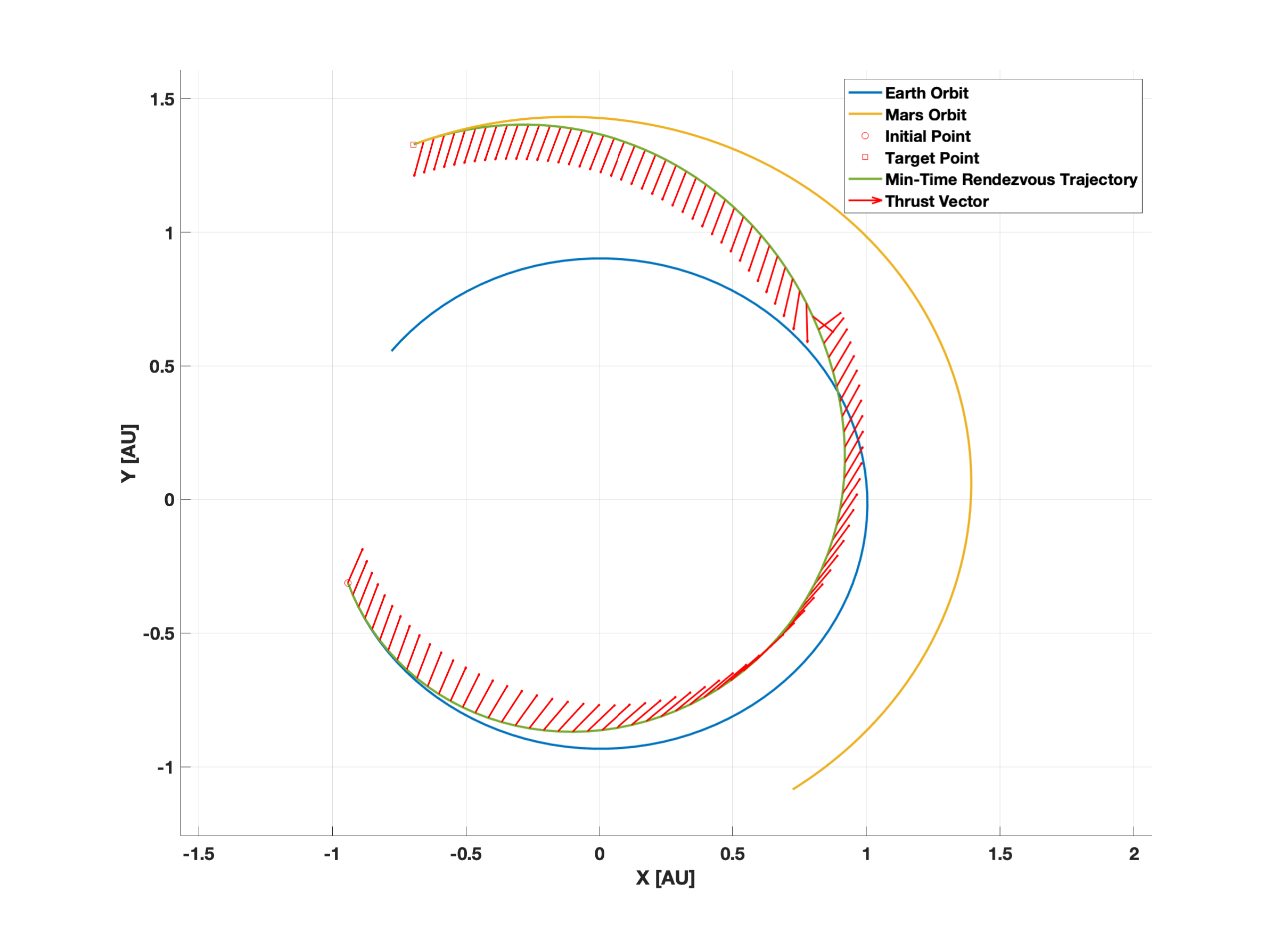}
  \label{fig:EMdirecttraj}
}

\subfloat[Componenets of the control vector.]{%
  \includegraphics[clip,width=0.9\columnwidth]{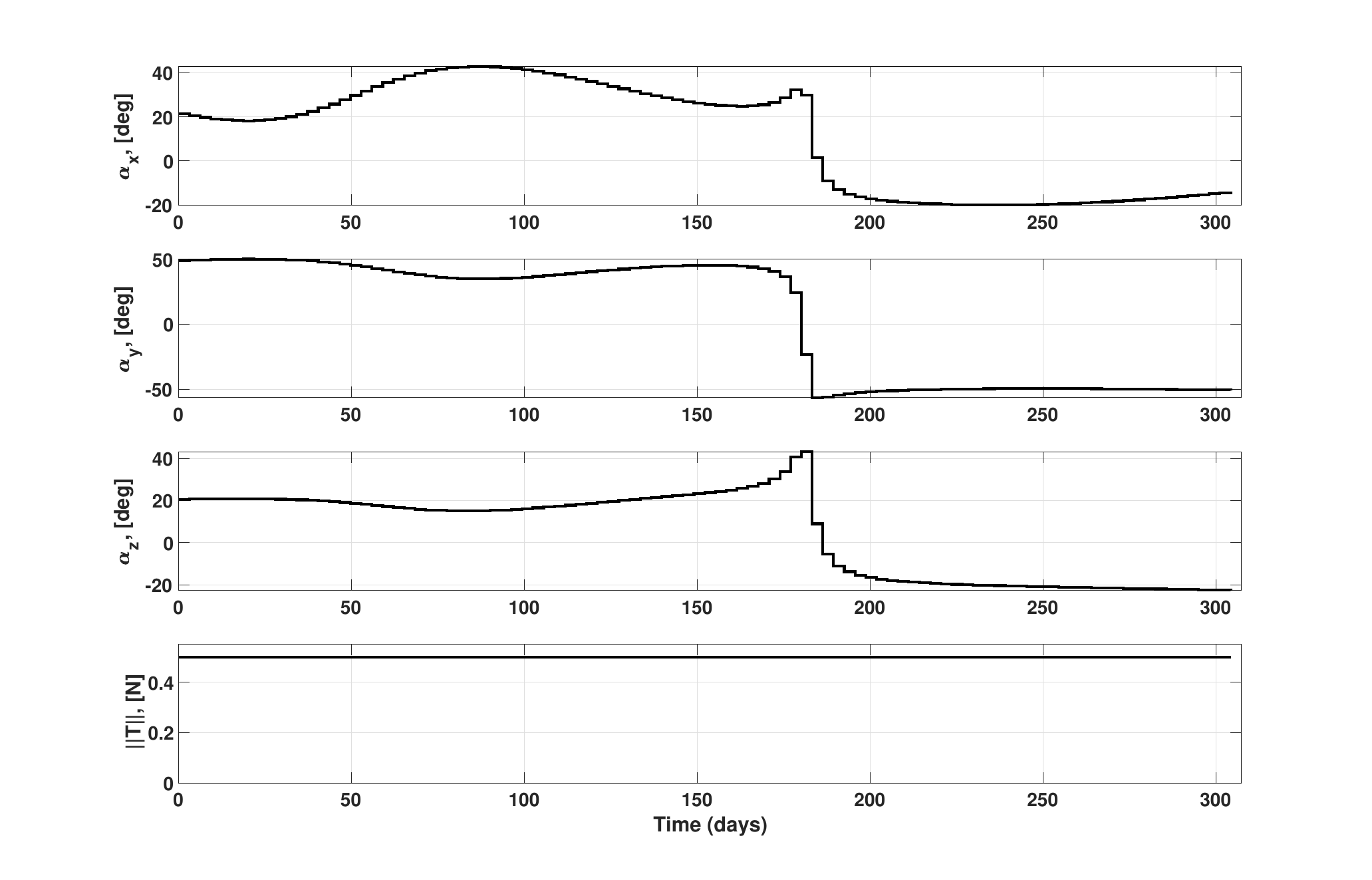}
  \label{fig:EMdirectcont}
}

\caption{Minimum-time Earth-Mars problem solved with CasaDi.}
\label{fig:EMdirectoptimal}
\end{figure}

We seek to validate the accuracy of the reachable set estimate by comparing the results in Sec.~\ref{sec:SRSA} with a  minimum-time trajectory optimization problem. We solve a minimum-time direct OCP using the same initial conditions listed above. The minimum-time Earth-Mars OCP was solved using CasADi \cite{andersson2019casadi}), the resulting optimal trajectory and control are presented in Figure \ref{fig:EMdirectoptimal}. Observe that the propulsion system is operating at maximum thrust for the duration of the flight, and only the thrust direction (red arrows) in Figure \ref{fig:EMdirecttraj} varies. The minimum-time Earth-Mars problem, when solved, resulted in a trajectory with a time of flight of 307 days. The trajectory is also leveraging Oberth effect (i.e., it falls into the gravitational well of the Sun to maximize the rate of the change of its energy). 

% Taheri posed a minimum-fuel Earth-Mars problem in \cite{taheri_enhanced_2016,taheri_performance_2018} and described how a (Hamiltonian) TPBVP arises from an indirect formulation. This TPBVP can be solved with a single-shooting method (e.g., using MATLAB's built-in \texttt{fsolve} function). We use this same formulation and manually iterate on time of flight until a trajectory is determined that produces a constant maximum thrust profile in the shortest time of flight. The same solution can also be obtained by formulating a minimum-time OCP and using an NLP solver (e.g. CasADi \cite{andersson2019casadi}). 

The solution of the considered minimum-time version of the Earth-Mars problem aligns with an emphemeris-consistent rendezvous boundary conditions (i.e., the position and velocity of Mars depend on the time of flight). The resulting minimum-time trajectory can be used to demonstrate the validity of the reachable set estimations because a minimum-time trajectory, by definition, must terminate on the reachable set. To show this, we compute the reachable position and velocity sets for the time of flight computed by the optimal single-trajectory indirect OCP \cite{taheri2017co}. Refer to Figure \ref{fig:velevo} that shows the evolution of the velocity reachable set through flight times from 120 to 180 days. We select these time of flights since they correspond to the thruster direction reversing in Figure \ref{fig:EMdirectoptimal}. The reachable set appears to fold upon itself when chronologically stepping from 140 days (Fig. \ref{fig:vel140}) to 160 days (Fig. \ref{fig:vel160}). This shows that the behaviors of the reachable position and velocity sets are not obviously correlated.
\begin{figure}[!tbp]
\centering
\subfloat[120 Days.]{{\includegraphics[clip, width=.475\columnwidth]{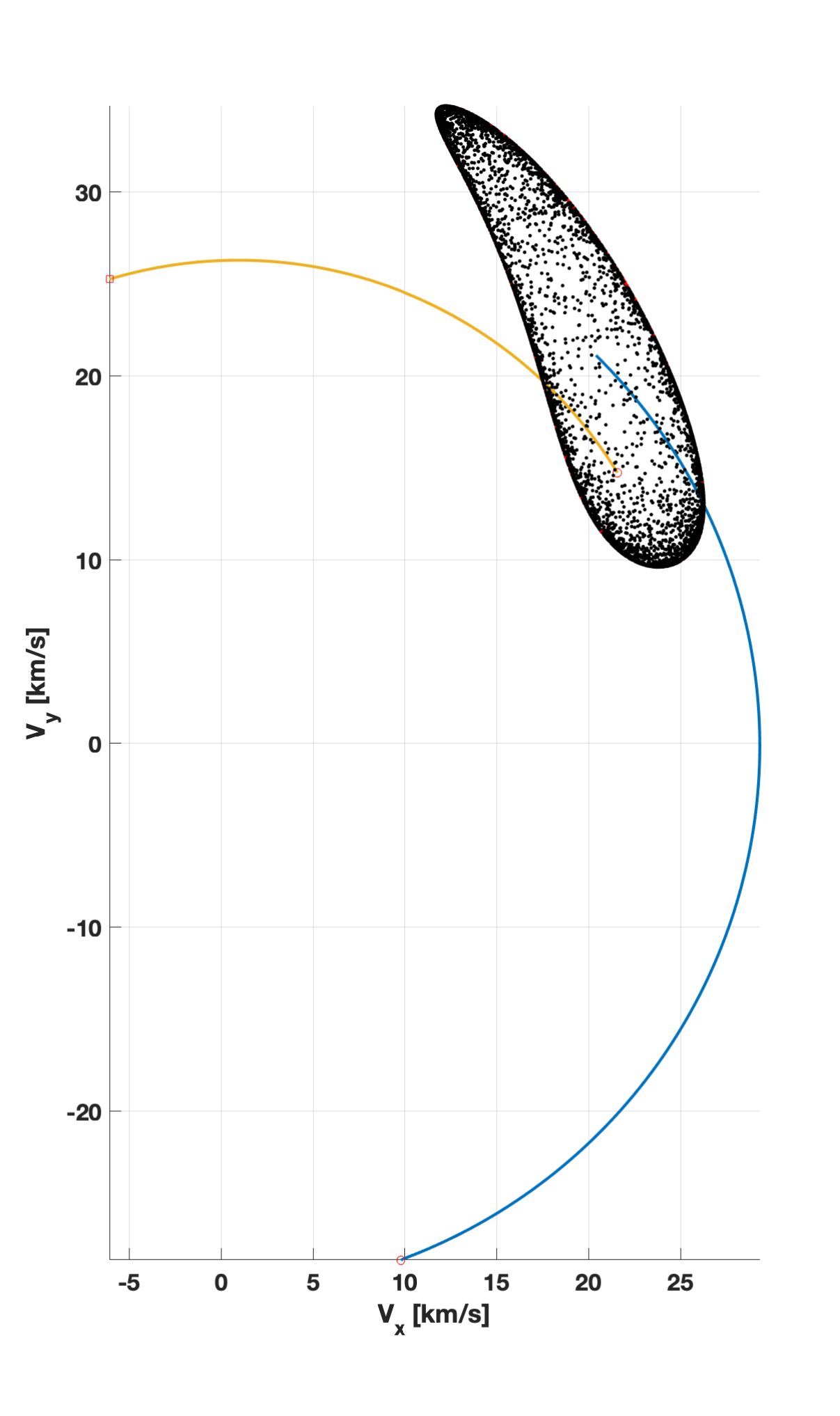}}\label{fig:vel120}
}  
\subfloat[140 Days.]{{\includegraphics[clip, width=.475\columnwidth]{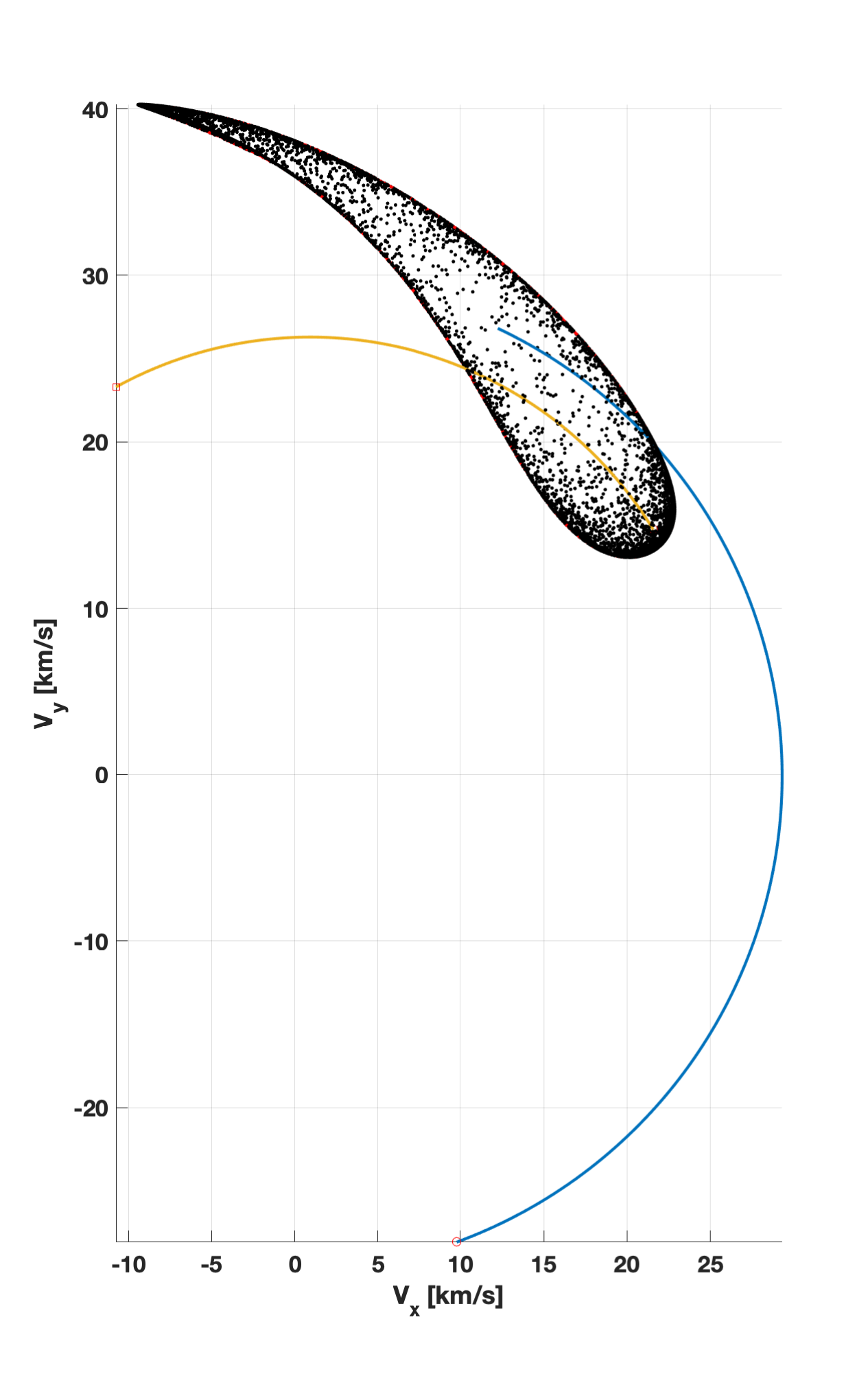}}\label{fig:vel140}
} 
\end{figure}

\begin{figure}[!tbp]
\centering
\ContinuedFloat
\subfloat[160 Days.]{{\includegraphics[clip, width=.7\columnwidth]{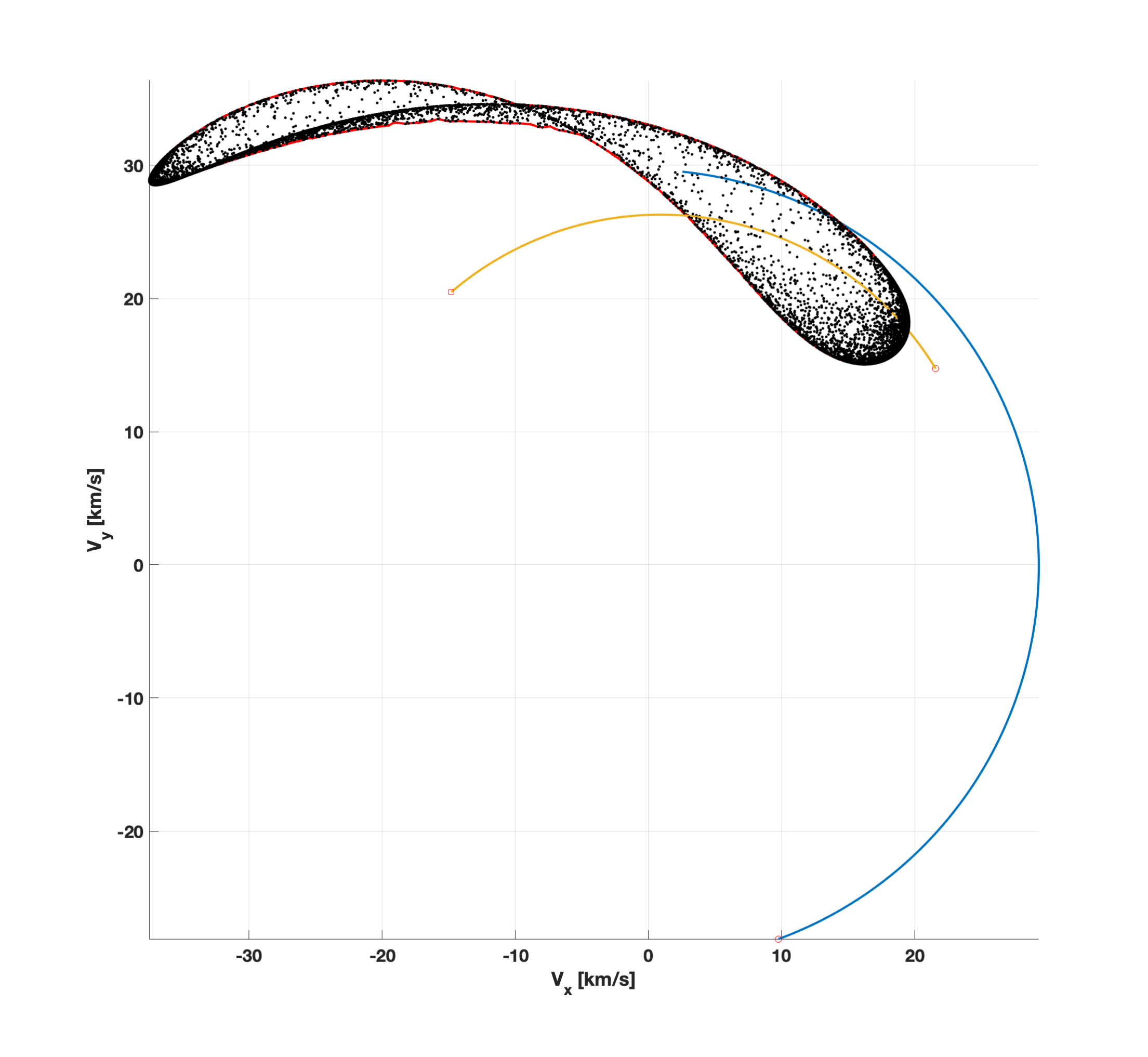}}\label{fig:vel160}
  }
  
\subfloat[180 Days.]{{\includegraphics[clip, width=.7\columnwidth]{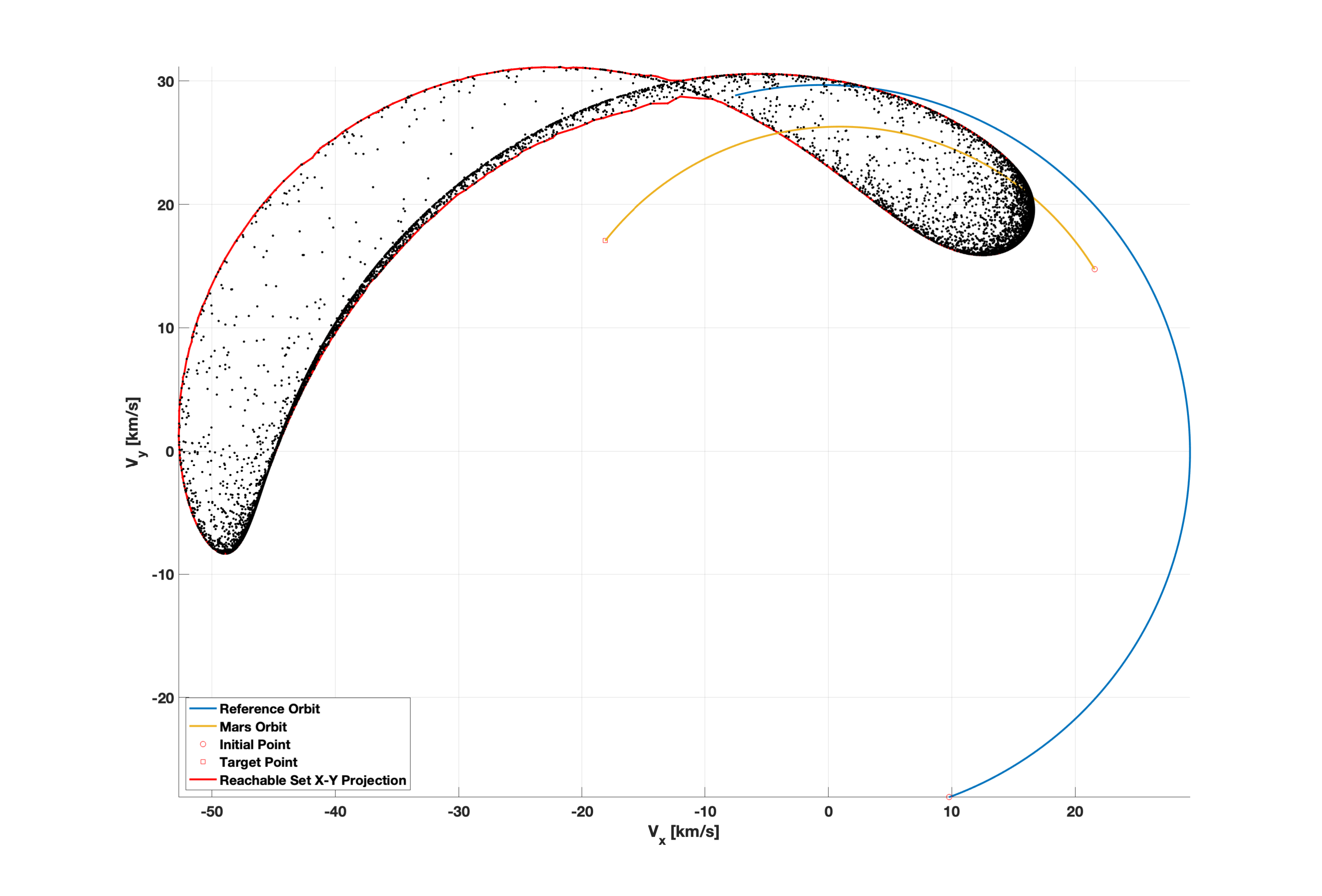}}\label{fig:vel180}
}

\caption{Evolution of reachable velocity set vs. different time horizons.}
\label{fig:velevo}
\end{figure}

\begin{figure}[!tbp]
\centering 
\subfloat[Position reachable set.]{%
  \includegraphics[clip,width=0.95\columnwidth]{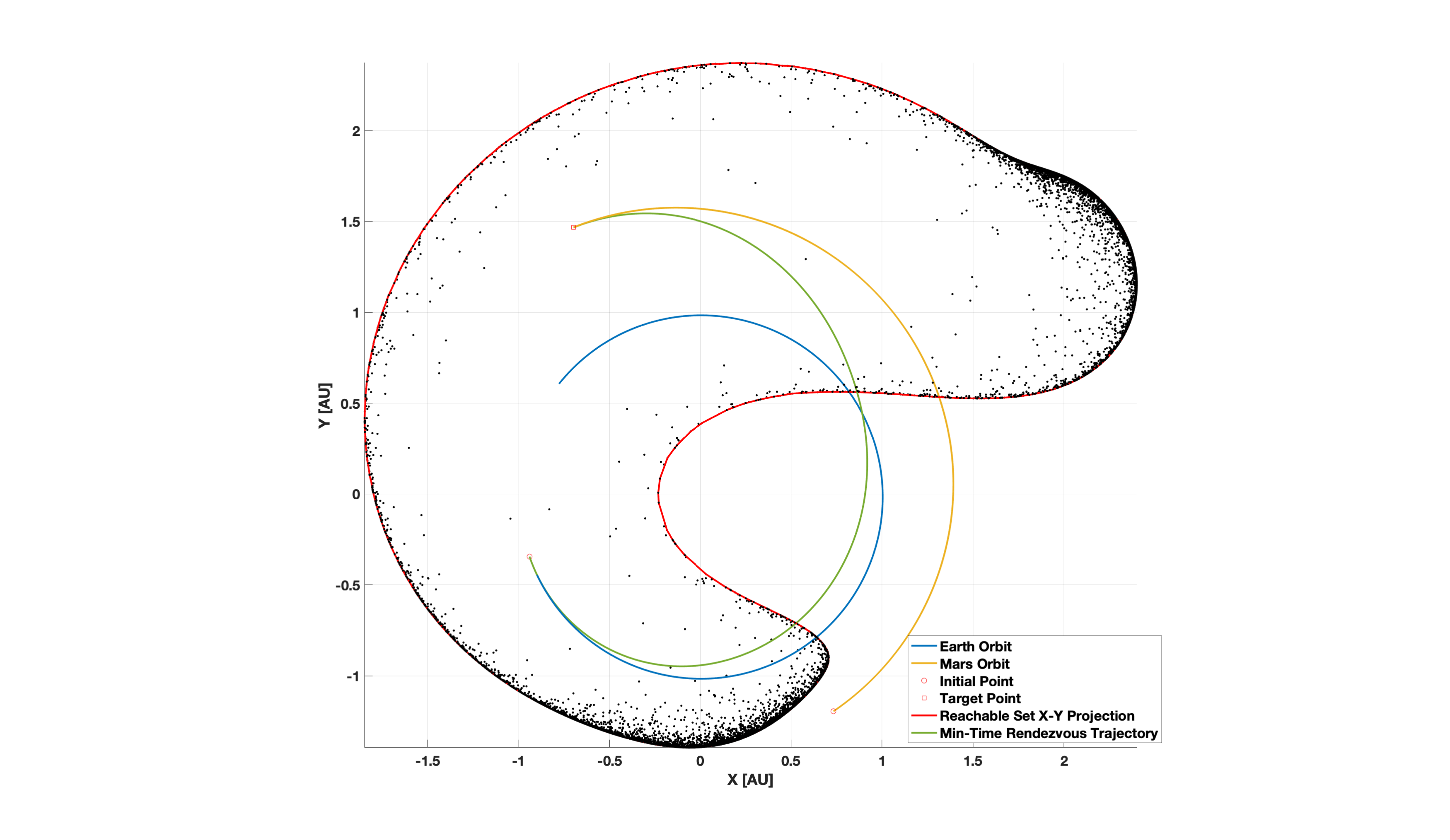}
  \label{fig:posreach}
}

\subfloat[Velocity reachable set.]{%
  \includegraphics[clip,width=0.95\columnwidth]{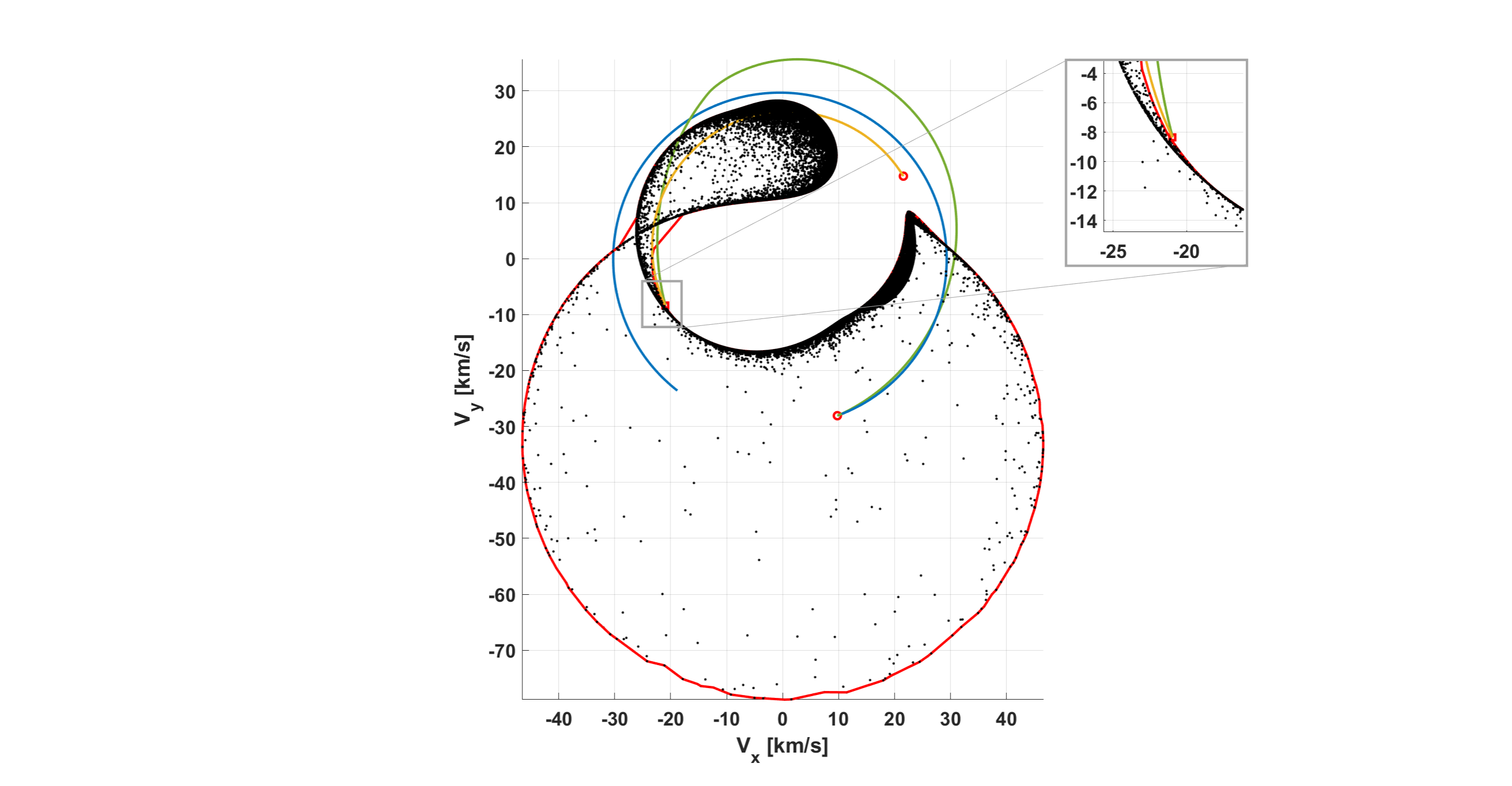}
  \label{fig:velreach}
}

\caption{Reachable set over a 307-day time horizon with 5000 sample trajectories and minimum-time trajectory.}
\label{fig:2bodyposvel}
\end{figure}

The reachable set algorithm was used to compute the reachable set for a time horizon of 307 days with 5000 sample points. Figure \ref{fig:posreach} shows the position reachable set for the given initial conditions. Upon examination, we see that the indirect OCP generates a position-feasible trajectory (green line) due to the fact that in position-space, the minimum-time trajectory terminates on the interior of the projected position reachable set. Figure \ref{fig:velreach} depicts the velocity reachable set. The indirect minimum-time  trajectory terminates on the boundary of the velocity reachable set. This indicates that for the considered parameters of the propulsion system the velocity states are more restrictive, compared to the position set, and reduces the combined reachable velocity. Moreover, the results are consistent and validates the accuracy of the method for determining an accurate estimate of reachable set for low-thrust spacecraft. However, the sensitivity of the reachable set for fixed-time rendezvous maneuvers is more nuanced, as is discussed for a rendezvous maneuver from Earth to asteroid 1989ML in \cite{taheri_how_2020}, which is a consequence of the deformed reachable set with the orbit of Mars.

\subsection{Reachable set propagation and analysis of the compuation time}
One useful application of this reachable set determination method is the ability to rapidly determine if a planet/object can be encountered from given initial conditions in a specified time horizon. To demonstrate this capability, the reachable set algorithm was initialized such that one stage is equal to one day. Then, 10,000 sample trajectories are computed for each trial. Data indicating the time of flight and simulation run times are presented in Table \ref{reach1}. The fast computation time allows for rapid iteration to ``brute force'' search for a time of flight that allows a target to become reachable.

% miss distance, or the closest distance to the target achieved by one of the sample trajectories, the reachable position set volume corresponding to a hull volume interpolated from sample points, 

\begin{table} [ht!]
\caption{\label{tab:table1} Earth-Mars reachable set propagation with 10,000 sample points.}
\centering
\begin{tabular}{|c|c|}
\hline
Time of Flight  & Computation time\\
(days) & (seconds)\\
\hline
100 & 9.0599\\
\hline
150 & 17.3945\\
\hline
200 & 18.1148\\
\hline
250 & 21.8414\\
\hline
300 & 26.2118\\
\hline
\end{tabular}
\label{reach1}
\end{table}

\begin{figure}[hbt!]
\centering
\includegraphics[scale = 0.5]{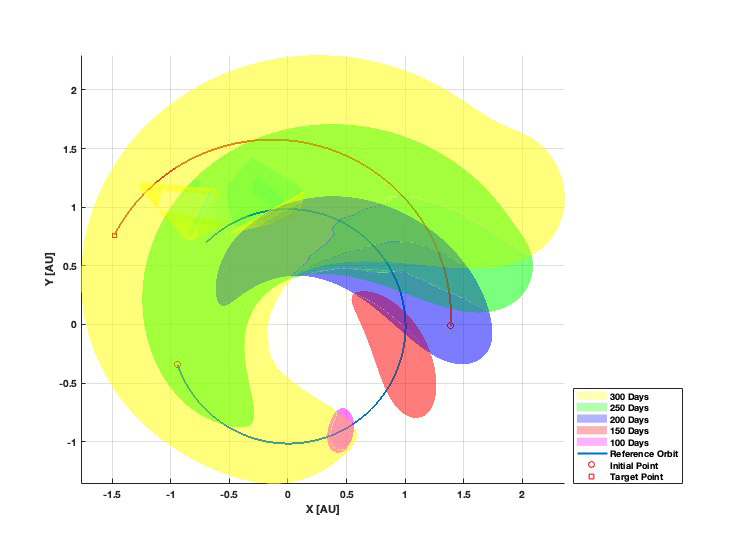}
\caption{Depiction of the Earth-Mars low-thrust position-level of the reachable set.}
\label{setprop}
\end{figure}

Figure \ref{setprop} indicates the expansion of the reachable position set through an increasing time horizon. This pseud-qualitative plot was generated by interpolating all sampled trajectories on the reachable set to graphically illustrate the reachable set. The reachable set algorithm was run five times, with flight times ranging from 100 to 300 days. The trajectory of Mars (red line) shows how the target position evolves with increasing time of flight. 
The power of the rapid reachable algorithm is illustrated here, since the successive generation of many reachable sets for complex problems is extremely computationally expensive. %The technique in this paper allows for rapid generation of these sets to quickly determine if an intercept occurs. 
From Figure \ref{setprop} we can see that time of flight must be greater than 250 days for Mars to be reachable. Recall that position-reachability is the first step in determining if an object can be reached if a rendezvous scenario is of interest, which confirms the results of the previous section. %The future work on this front will be to add an outer loop onto the reachable set algorithm to iteratively determine the minimum time of flight required for intercepts from given initial conditions. From this data we will be able to offer comparisons between traditional optimal control formulations with this reachability approach.
% Figure \ref{volume} plots data from Table \ref{reach1} to indicate the exponential growth rate of the reachable position set.

% \begin{figure}[hbt!]
% \centering
% \includegraphics[scale = .5]{Figures/volumeplot.pdf}
% \caption{Position reachable set volume vs. time for the Earth-Mars low-thrust problem.}
% \label{volume}
% \end{figure}

\subsection{Reachable set variations with uncertainties in initial boundary conditions}
We revisit the simulation presented in Section \ref{sec:SRSA} to demonstrate the ability to add uncertainties in initial boundary conditions and to assess the impact of them on the reachable set. Two types of initial boundary conditions are considered: 1) position and velocity uncertainty, and 2) an impulse maneuver. Both types of constraints can be incorporated simultaneously, however, we separate them in these results to demonstrate their individual effects. We consider departure from the Earth (reference) orbit.

\begin{figure}[!tbp]
\centering 
\subfloat[With uncertainties on the initial position and velocity vectors (per defintion given in Eq.~\eqref{posellip}).]{%
  \includegraphics[clip,width=0.95\columnwidth]{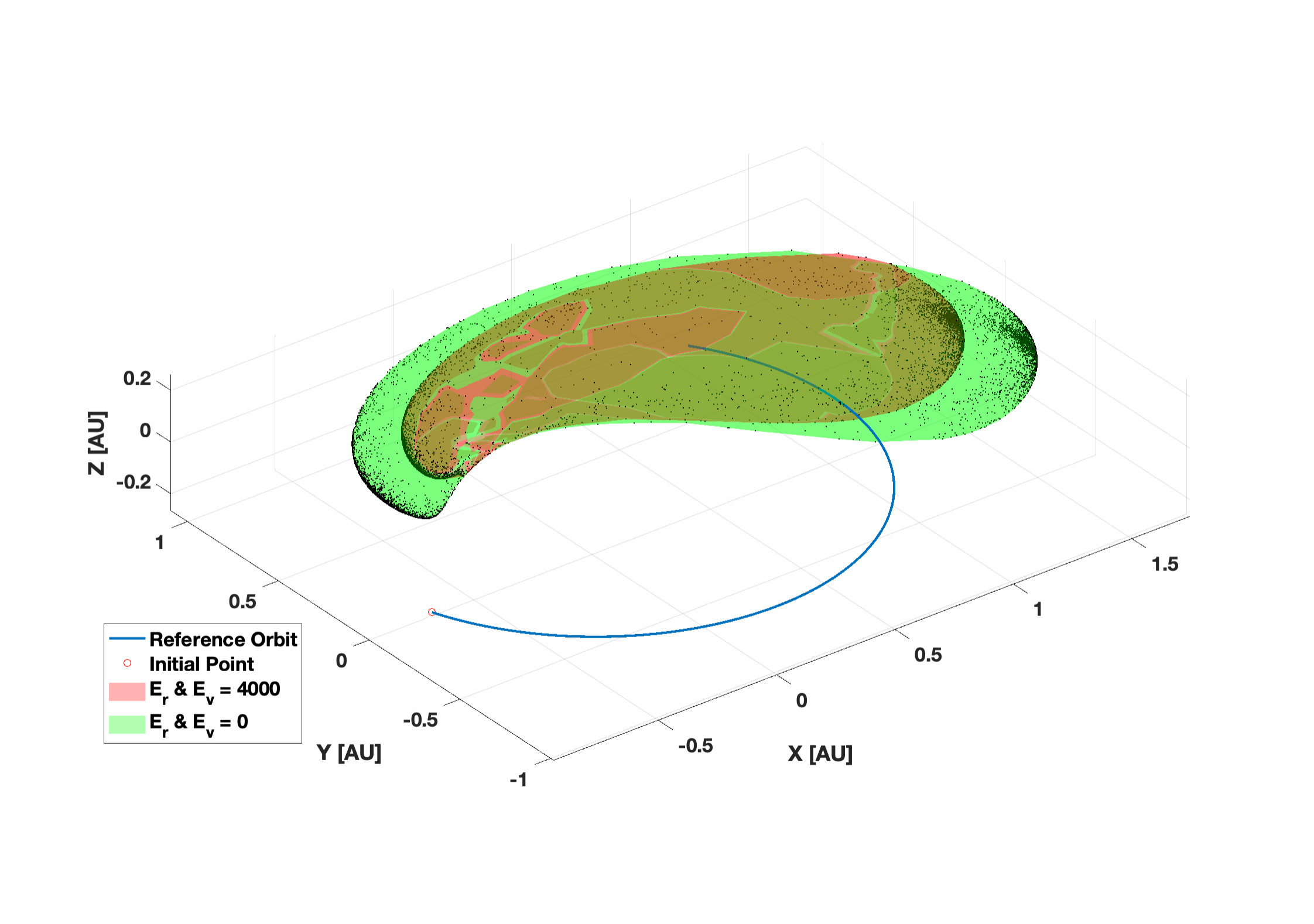}
  \label{fig:posvelu}
}

\subfloat[With an initial impulse maneuver (per defintions given in Eqs.~\eqref{dv2bind} and \eqref{dv2nobind}).]{%
  \includegraphics[clip,width=0.95\columnwidth]{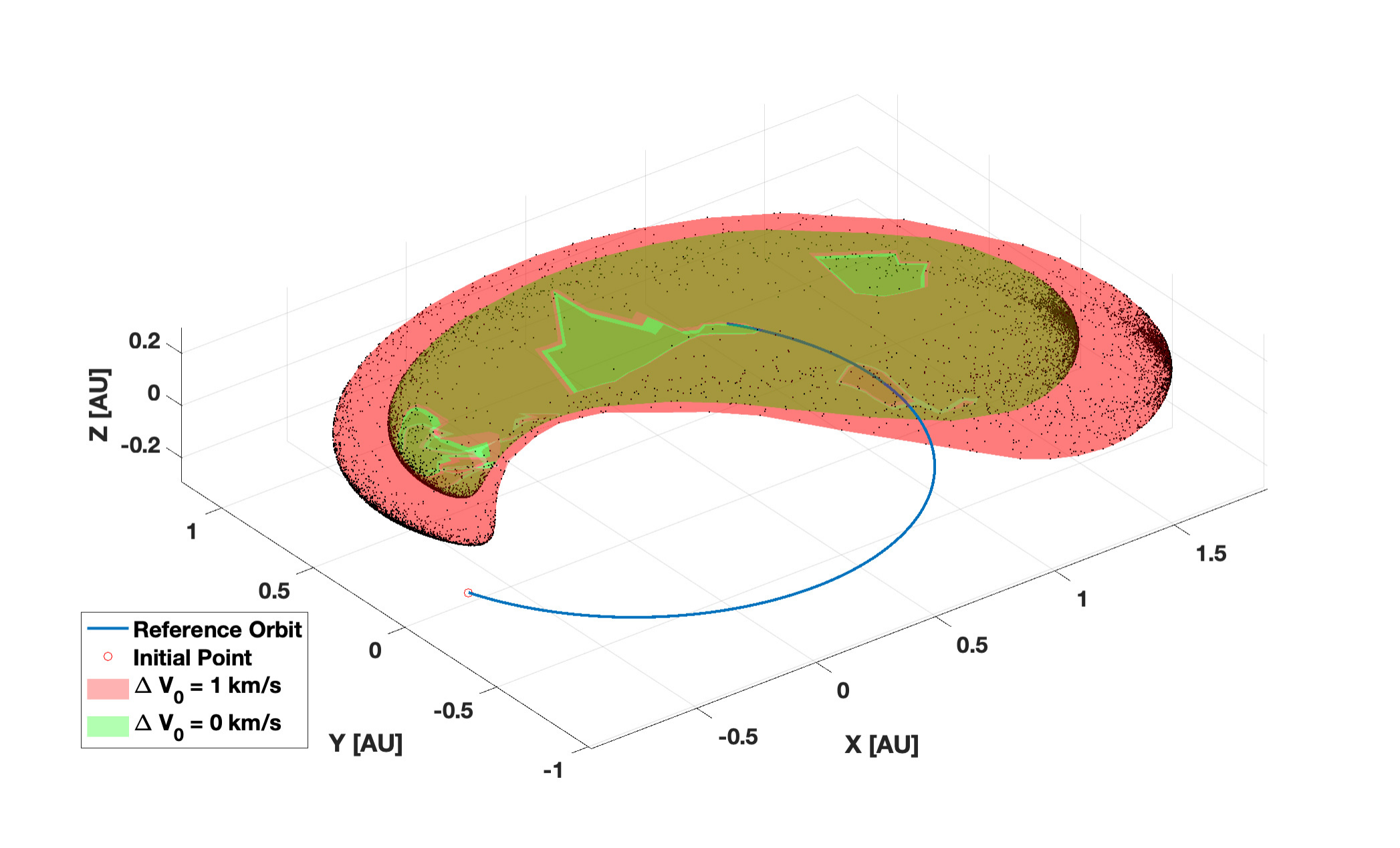}
  \label{fig:impulse}
}

\caption{Reachable sets over a 200-day time horizon with 5000 samples  with different initial types of uncertainties.}
\label{fig:2bodyuncertainty}
\end{figure}

Figure \ref{fig:posvelu} shows the reachable set for a case without position and velocity uncertainty (green surface) and a case for ellipsoidal position and velocity (red surface). The addition of ellipsoidal uncertainty constraints clearly affects the nominal reachable set, in that it appears the red set had decreased in size relative to the green set. However, this is not necessarily the case since the uncertainty constraint has also shifted the red reachable set out of plane relative to the green reachable set, making it appear smaller. The main conclusion is the demonstration that ellipsoidal uncertainty constraints can significantly affect the reachable set. Figure \ref{fig:impulse} presents the expected result of adding an initial $\Delta V$ impulse (which can be considered as excess velocity provided by a launch vehicle). The impulse maneuver excites the spacecraft to a higher energy state, which ultimately leads to a larger reachable set. The red set of the impulse maneuver, is noticeably expanded after a 200-day time interval compared to the zero-impulse maneuver reachable set (green color). 

\section{CR3BP Reachable Set Results} \label{sec:cr3bpresults}

%Results have been obtained for a low-thrust spacecraft in cislunar region under the CR3BP dynamics. 
Results are presented for three cases by considering different spacecraft parameters over various time horizon values to showcase the versatility of the proposed algorithm and to reveal the characteristics of the resulting reachable sets. Case 1 is an L1 point initial condition with no velocity over increasing time horizon values and is in agreement with results presented in \cite{jones_reachable_2023}. Case 2 simulates an L2 Halo orbit over both short- and long-time-horizon values and illustrates the chaotic natural dynamics and reachable set evolution of the CR3BP. Case 3 updates the initial condition to a 9:2 NRHO, which is the chosen orbit for NASA's Lunar Gateway station \cite{trofimov2020transfers}. Case 1 serves to validate the reachable set results, whereas Cases 2 and 3 present novel reachable set solutions in cislunar space. Last, a theoretical explanation involving invariant manifolds is presented to offer additional insights into the evolution of the reachable set in the CR3BP. In all results presented, the Earth and Moon are not illustrated to scale and the segment duration, $\Delta t$, was determined empirically as the time of flight normalized by characteristic time divided by 200.

\subsection{L1 Point reachability}

\begin{figure}[!tbp]
\centering 
\subfloat[50-100 Hour Time of Flight]{%
  \includegraphics[clip,width=0.8\columnwidth]{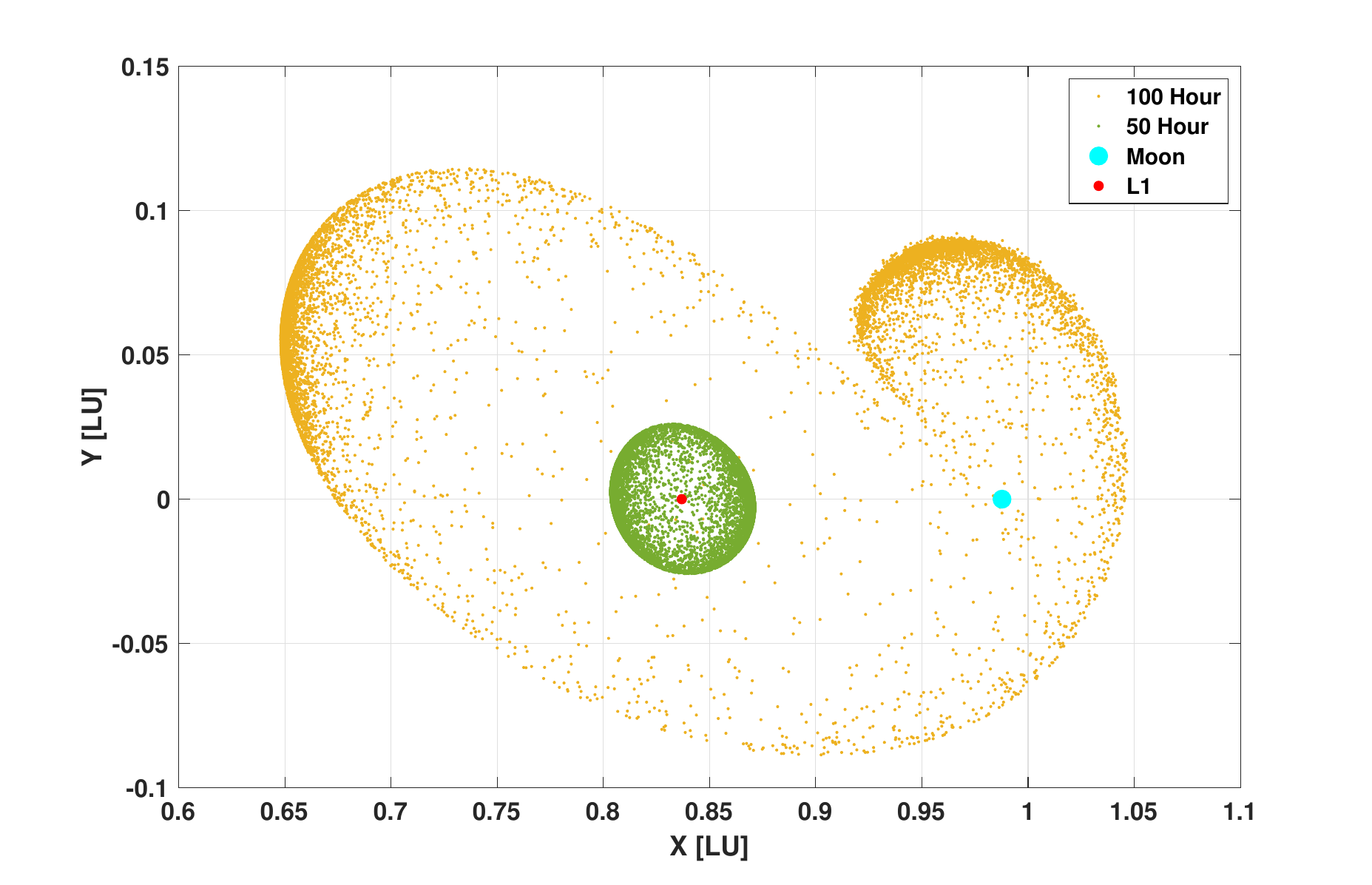}
  \label{fig:L1short}
  }

\subfloat[150-200 hour time of flight.]{%
  \includegraphics[clip,width=0.8\columnwidth]{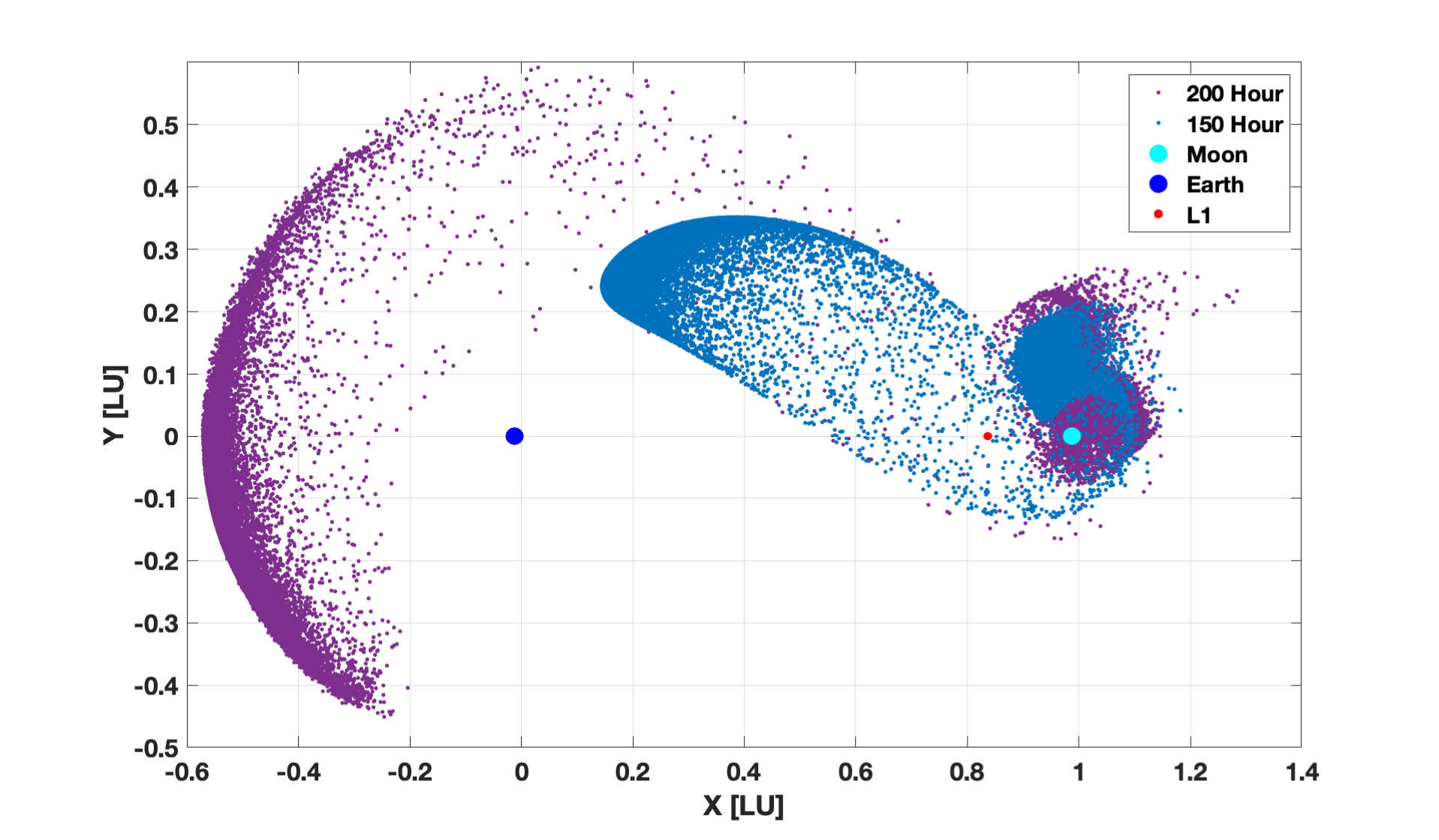}
  \label{fig:L1long}
  }

\caption{L1 Point reachable sets with increasing time of flight.}
\label{fig:L1main}
\end{figure}

A low-thrust spacecraft with $T_\text{max} = 1$ Newton, $I_\text{sp} = 2000$ seconds, and 1500 kg of initial mass is considered. The initial state vector is given in nondimensionalized coordinates as  $\boldsymbol x = [0.836892919; 0; 0; 0; 0; 0]^{\top}$. The reachable set was computed for times of flight ranging from 50 to 200 hours. Each trial uses 100,000 sample trajectories. The longest simulation run time, for 200 hours, took 119 seconds. Comparing with \cite{jones_reachable_2023}, which is a planar implementation of cislunar CR3BP, Figure \ref{fig:L1main} utilizes full three-dimensional dynamics. The sample points that lie on the interior of the reachable set with this $X-Y$ projection are due to the out-of-plane effects and are simply reachable points with a $Z$-coordinate not equal to zero. The results validate the accuracy of the linearization of the CR3BP dynamics. 

% \begin{figure}[hbt!]
% \centering
% \includegraphics[scale = .34]{Figures/L1dtf.pdf}
% \caption{L1 Point Reachable Sets with Increasing Time of Flight.}
% \label{L1main}
% \end{figure}

 \subsection{L2 Halo orbit reachability}

 Case 2 considers computation of the reachable set for a L2 Halo Orbit. A low-thrust spacecraft with $T_\text{max} = 0.2$ Newton, $I_\text{sp} = 3000$ seconds, and $m_0 = 1000$ kg is considered. The initial state vector is given in nondimensionalized coordinates as $\boldsymbol x = [1.17204419281306; 0; -0.0862093101977581; 0; -0.188009087163036; 0]^{\top}$. The period of the L2 reference orbit is $346.322857$ hours. The 3D trajectory  encounters significant multi-body effects. The reachable set was computed for times of flight of 150 and 350 hours. The 150-hour time horizon demonstrates the evolution of the reachable set in the vicinity of the moon, and the 350-hour trial shows growth and a disjoint reachable set.

 \begin{figure}[!tbp]
\centering 
\subfloat[3D view.]{%
  \includegraphics[clip,width=0.8\columnwidth]{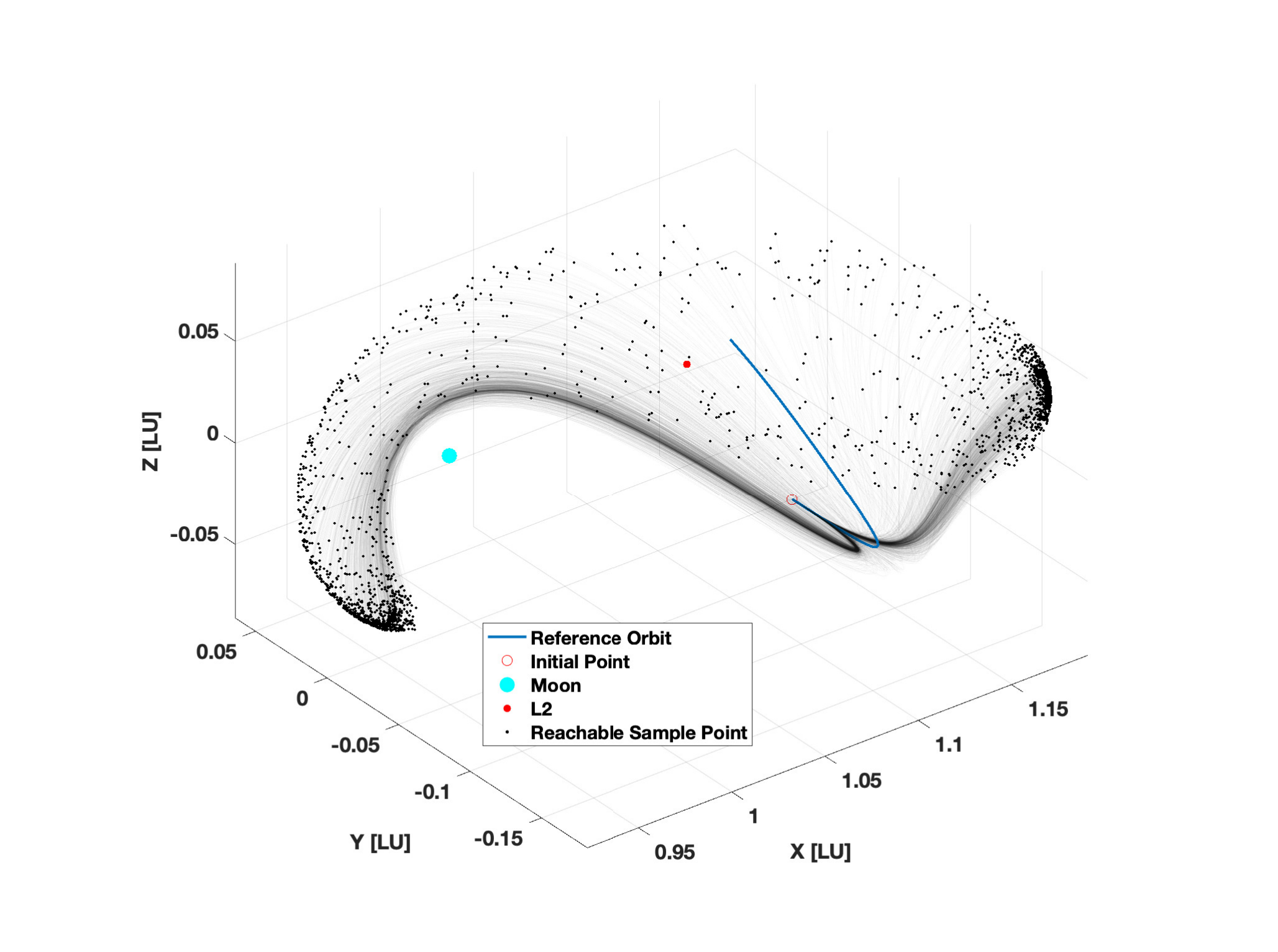}
  \label{fig:L2xyz150}
  }

\subfloat[XY view.]{%
  \includegraphics[clip,width=0.8\columnwidth]{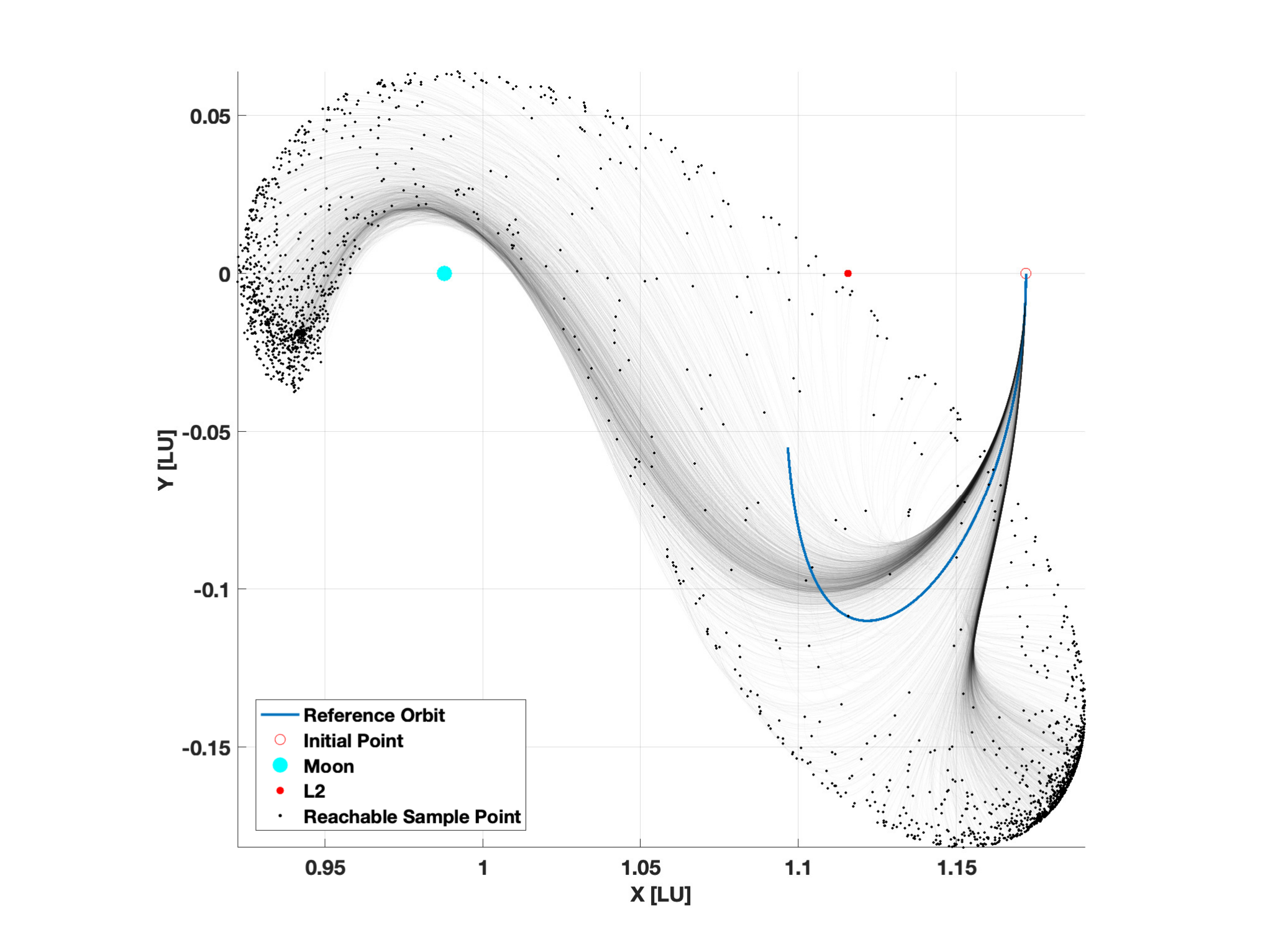}
  \label{fig:L2xy150}
  }

% \subfloat[XZ view]{%
%   \includegraphics[clip,width=0.95\columnwidth]{Figures/L2_xz_150h_2e3.pdf}
%   \label{fig:L2xz150}
%   }

% \subfloat[YZ view]{%
%   \includegraphics[clip,width=0.95\columnwidth]{Figures/L2_yz_150h_2e3.pdf}
%   \label{fig:L2yz150}
%   }
\caption{Position reachable set over a 150-hour (6.25 days) time horizon with 2000 sample trajectories.}
\label{fig:L2150h}
\end{figure}

Figure \ref{fig:L2150h} shows the 150 hours expansion of the reachable set for a low-thrust spacecraft under CR3BP dynamics with 2,000 sample trajectories. Even on this short time horizon, one can identify that small changes in the thrust vector (especially during initial stages) directly contribute to a reachable set that is fully integrated in three dimensions. Unlike the two-body  results (see Figs. \ref{fig:EMreachabilitymain}, \ref{fig:2bodyposvel}, and \ref{fig:2bodyuncertainty}) in which it is straightforward to interpolate the reachable sample points to generate a reachable set boundary projection, it is very difficult to accomplish with the abstract reachable set shape that is illustrated here. Additionally, Figure \ref{fig:L2150h} shows an incipient tendency for the reachable trajectories to cluster into two bands, each propagating in opposite directions from the initial point. This phenomenon is the primary purpose of including the relatively long 350-hour time horizon case and will be explained in the section on invariant manifolds. 

\begin{figure}[!tbp]
\centering 
\subfloat[3D view.]{%
  \includegraphics[clip,width=0.8\columnwidth]{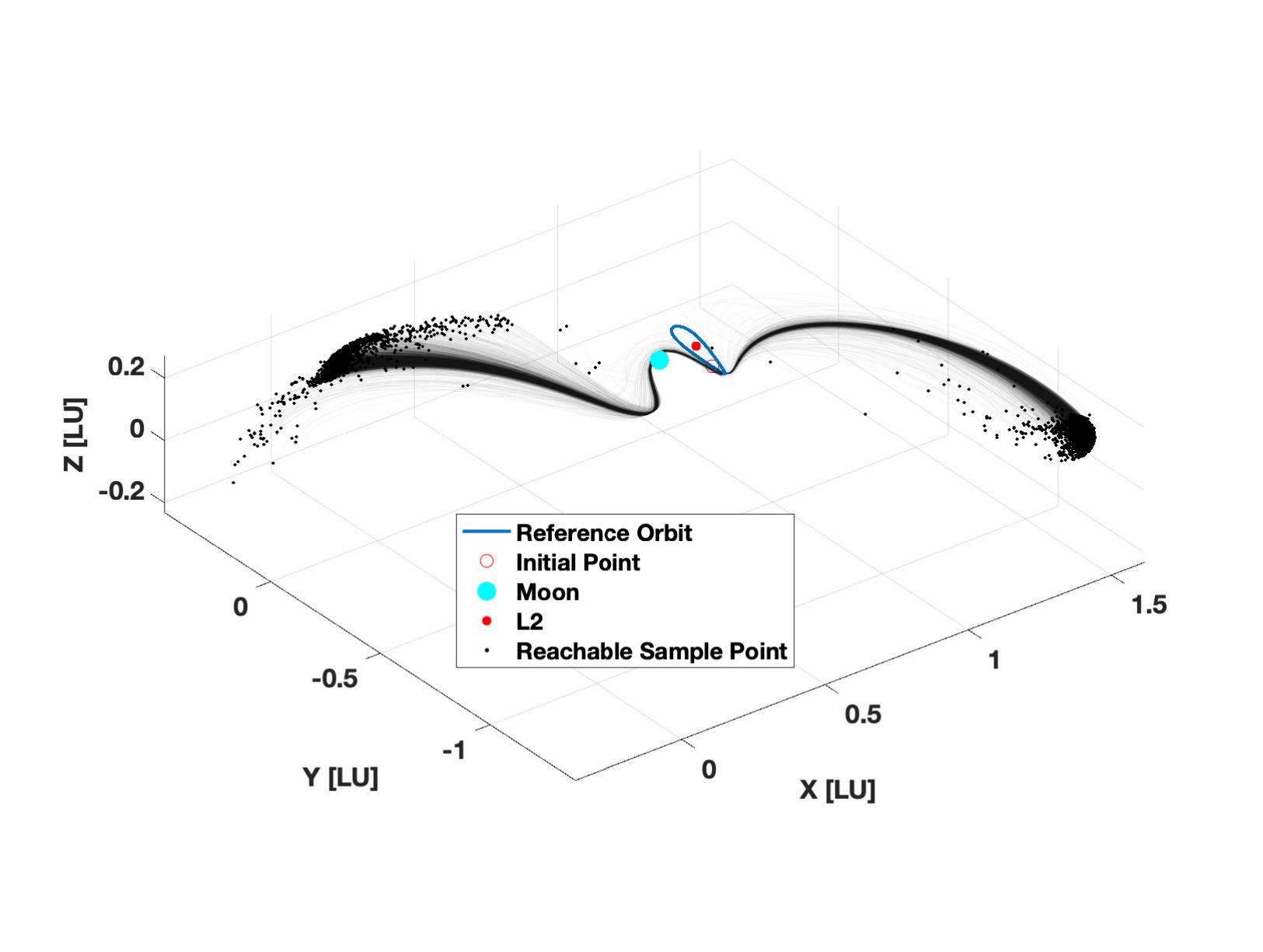}
  \label{fig:L2xyz350}
  }

\subfloat[XY view.]{%
  \includegraphics[clip,width=0.8\columnwidth]{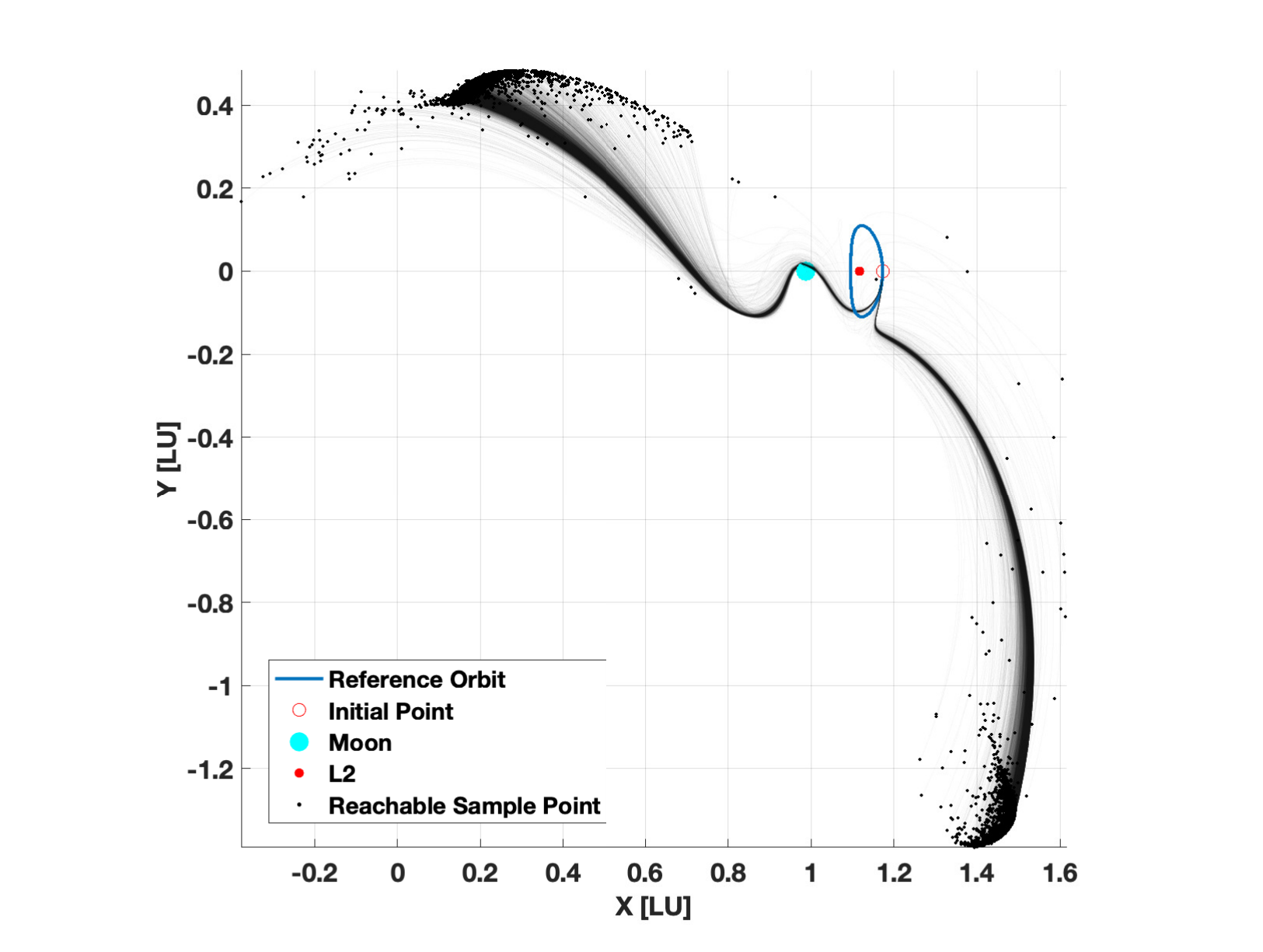}
  \label{fig:L2xy350}
  }

\caption{Position reachable set over a 350-hour time horizon with 10,000 sample trajectories.}
\label{fig:L2350h}
\end{figure}

For the 350-hour test case, the number of sample points was increased to 10,000 to increase the resolution of the captured reachable set. Figure \ref{fig:L2350h} presents the computed reachable set for the L2 Halo reference orbit. These results show clearly that the reachable set has disconnected into two primary clusters. One cluster completes a close fly-by with the moon and the second cluster exits the near vicinity of the moon into an area dominated by the gravitational attraction of the Earth. The section detailing invariant manifolds will outline the theoretical foundation for why this reachable set is disjoint and why the uniformly sampled trajectories are exhibiting unexpected clustering behavior.

\subsection{9:2 NRHO reachability}

Case 3 considers the reachable set for the Lunar Gateway earmarked 9:2 NRHO. A low-thrust spacecraft with $T_\text{max} = 0.2$ Newton, $I_\text{sp} = 3000$ seconds and $m_0 = 1000$ kg is used. The initial state vector is given in nondimensionalized coordinates as $\boldsymbol x = [1.0221, 0, -0.1821, 0, -0.1033, 0]^{\top}$\cite{thangavelu_transfers_2019}. The period of the L2 reference orbit is $157.500622$ hours. All results are computed with a 75-hour time of flight and 2,000 sample trajectories. 
%The 9:2 NRHO was chosen for this experiment since it has been chosen as the orbit for the NASA Gateway station. 
The region in the vicinity of the 9:2 NRHO will become a congested area in the near future, so it is imperative to be able to compute reachable sets to perform accurate spacecraft proximity operations in this space and execute potential collision-avoidance maneuvers.

\begin{figure}[!tbp]
  \centering
  \subfloat[3D view.]{\includegraphics[width=0.65\textwidth]{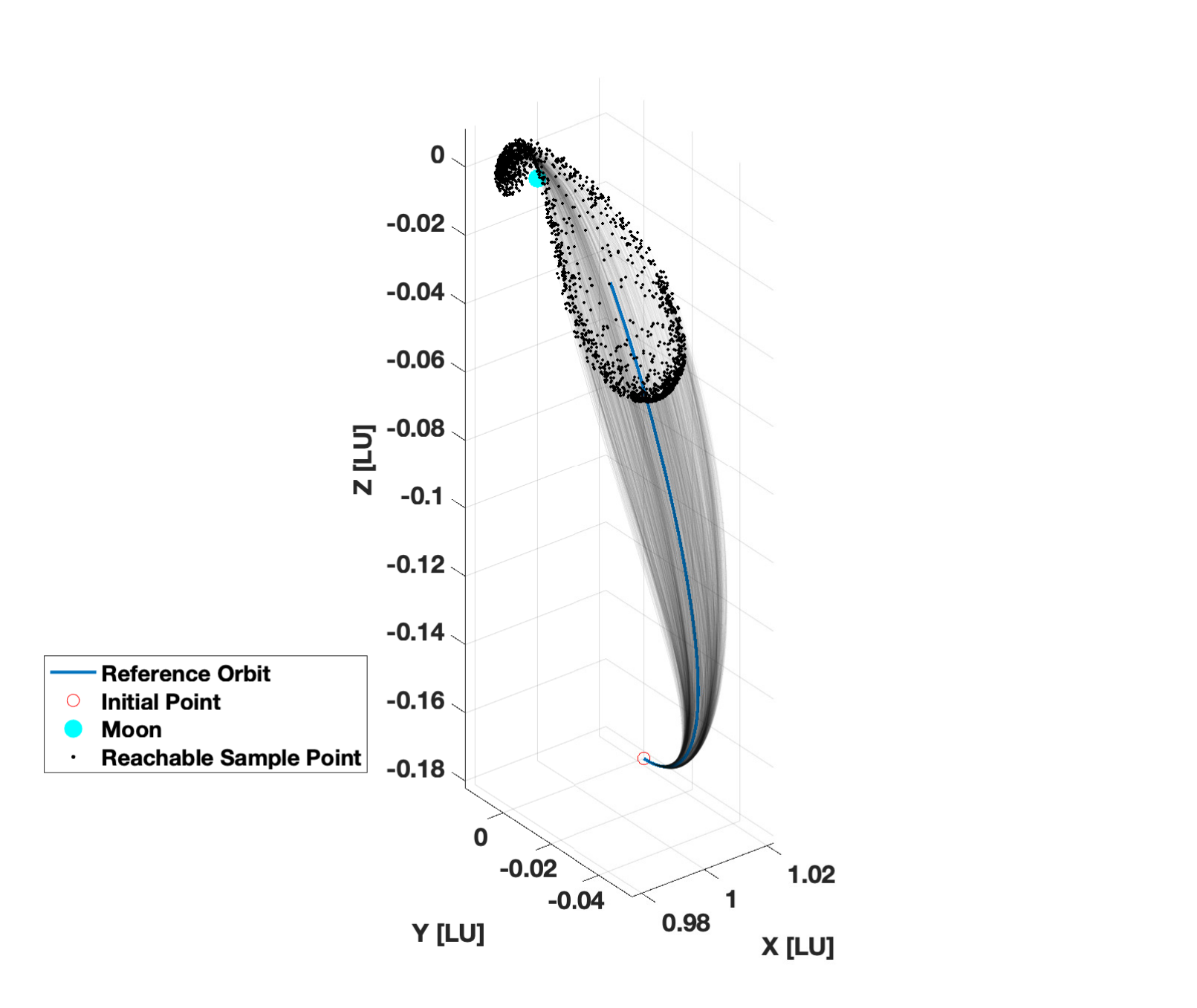}\label{fig:NRHOxyz}}
  \hfill
  \subfloat[XY view.]{\includegraphics[width=0.65\textwidth]{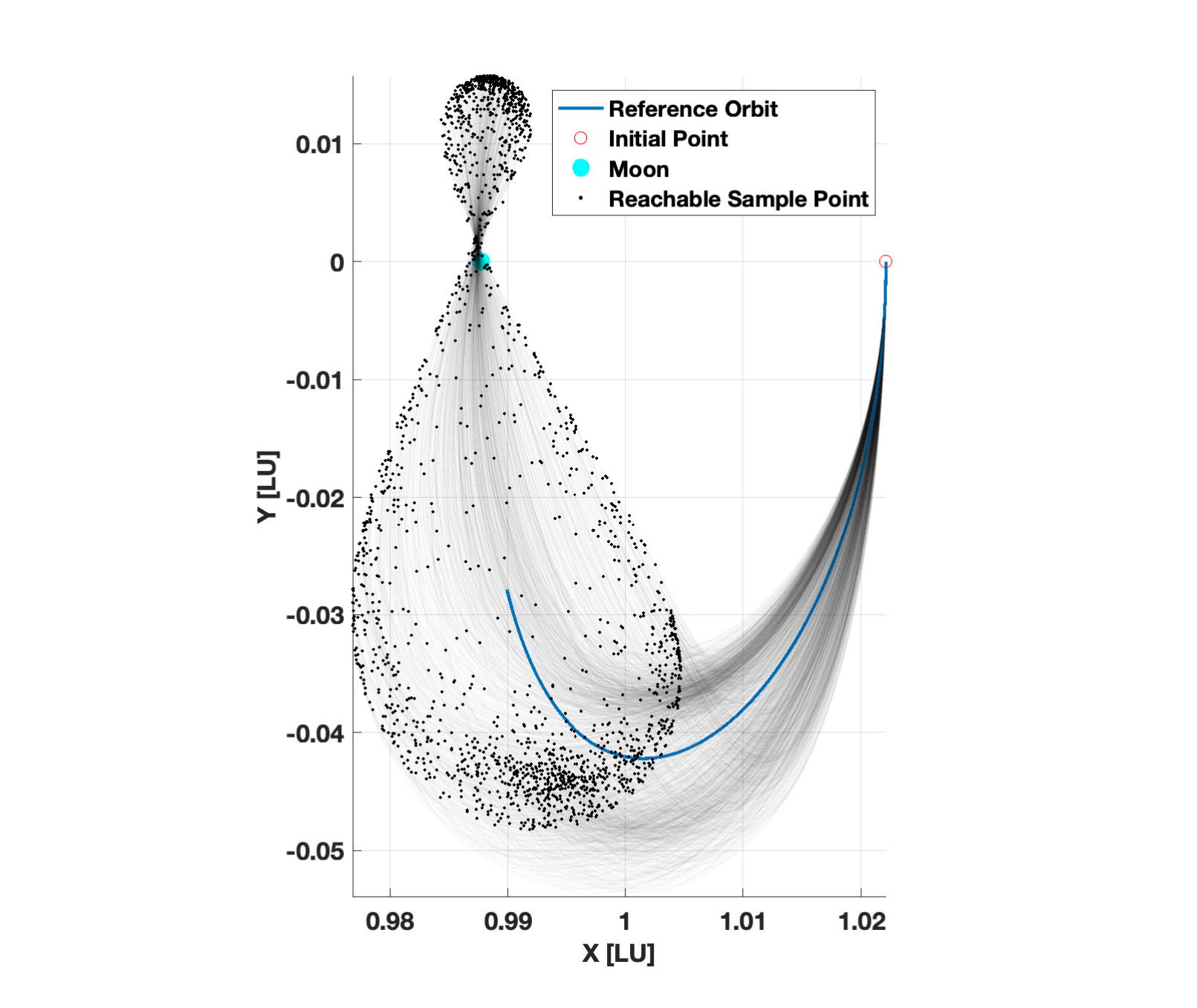} \label{fig:NRHOxy}}
  \caption{Position reachable set over a 75-hour time horizon with 2,000 sample trajectories.}
  \label{fig:NRHO}
\end{figure}

Figure \ref{fig:NRHO} presents the computed reachable set for the 9:2 NRHO over 75 hour time of flight with 2,000 sample trajectories. This result shows a new characteristic of the reachable set, that is, the reachable set not only expands, but also contracts with this time horizon. Before the spacecraft reaches the perilune of the 9:2 NRHO reference orbit, the reachable set expands in a manner consistent with previous results. At perilune, the natural CR3BP dynamics contract the reachable set before expanding again after perilune.

\subsection{Connection between reachable sets and invariant manifolds in CR3BP} \label{sec:invariantmanifolds}

Historically, low-energy transfers in three-body dynamical models originated with the idea of invariant manifolds. These dynamical structures exist in the vicinity of the three-body libation points and can be used to connect trajectories across vast distances in space \cite{mingotti2014combined}. An observation from the results presented for a low-thrust spacecraft in the CR3BP exhibits an interesting behavior over long time horizons, namely, the reachable set exhibits a bifurcation (i.e. there exist ``two'' different directions leading to reachable sets starting from the same initial condition). Additionally, when all the sampled minimum-time trajectories are plotted, they appear to be clustering around some unknown dynamical structure in the CR3BP. Examining the structure of the presented algorithm, the terminal costates are sampled from a uniformly distributed ball, which follows that the reconstructed thrust vector should also reflect the uniform sampling. This is not consistent with the observations, in that there appear some thrust vectors that are, broadly, ``more likely'' due to the chaotic nature of the  CR3BP. We will now attempt to explain this phenomenon of a disjoint reachable set and the clustering of sampled trajectories.

We briefly digress to a discussion from Chaos Theory and the bifurcation of the logistic map. The original 1967 paper from May \cite{may_simple_1976} details how chaos emerges from the iterative solution to the logistic difference equation when the growth rate exceeds 3.570. In detailed plots of the logistic map, examining the regions beyond 3.570 with an unstable fixed point, there appear to be some ``more'' stable points that appear as peaks of a probability distribution in the chaotic region. These are referred to as supertracks (ST) \cite{oblow_supertracks_1988}. Returning to the analysis of the reachable sets in the CR3BP, we hypothesize that the clustering of sampled reachable trajectories is analogous to ST's of the logistic map. The CR3BP is a chaotic system similar to the logistic map in that there exist stable manifolds and libation points. We ascertain that the clustering of sample trajectories is coincident with the invariant manifolds that exist as quasi-stable/unstable structures in the CR3BP. To empirically validate our hypothesis, invariant manifolds associated with the considered L2 Halo orbit will be computed. These invariant manifolds will be overlaid on the results from the reachable set analysis.

\begin{figure}[hbt!]
\centering
\includegraphics[scale = .27]{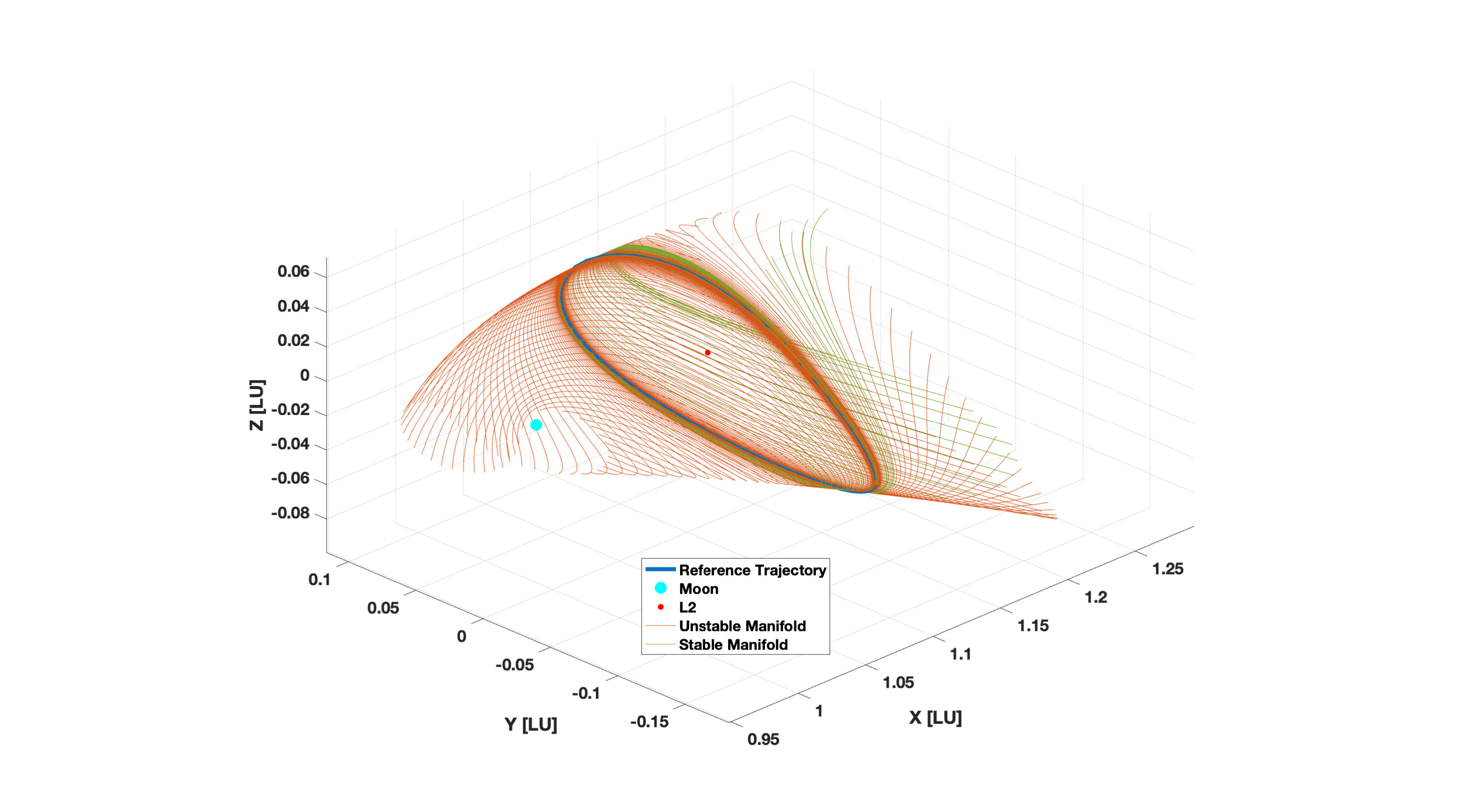}
\caption{Stable and unstable invariant manifolds of the Earth-Moon L2 Halo Orbit.}
\label{manifolds}
\end{figure}

Invariant manifolds are computed according to the process outlined in \cite{koon_dynamical_2011} and \cite{singh_eclipse-conscious_2021} and is straightforward to construct by the IMF. Beginning with the reference periodic orbit, we can compute the Monodromy matrix, $M(t_0)$, through a numerical integration of the state transition matrix, $\Phi$ (using Eq.~\eqref{cr3bp} with $\boldsymbol{u} = \textbf{0}$) from an initial point for one orbital period, $T$, as $M(t_0) = \Phi(t_0+T, t_0)$. Different selections of the initial point on the periodic orbit will yield different $M(t_0)$, which is useful for reconstructing the complete invariant manifolds. Next, compute the eigenvalues of $M(t_0)$. For Halo orbits, the six eigenvalues are $\lambda_1 = 1/\lambda_2$, $\lambda_3 = \lambda_4 = 1$, $\lambda_5 = \lambda_6^*$, where $\lambda_1$ and $\lambda_2$ are the real unstable and stable eigenvalues, respectively. In fact, $\lambda_3$ and $\lambda_4$ indicate neutral stability and $\lambda_5$ and $\lambda_6^*$ are a complex conjugate pair that lie on the unit circle and contribute only to rotation (i.e. neutrally stable directions). Extracting the eigenvectors associated with the unstable and stable eigenvalues allows for a perturbation in the unstable or stable direction on the initial condition as, $\boldsymbol{x_U}(t_0) = \boldsymbol{x}(t_0) + s \epsilon \boldsymbol{\xi_{U,0}}$, and $\boldsymbol{x_S}(t_0) = \boldsymbol{x}(t_0) + s\epsilon\boldsymbol{\xi_{S,0}}$
%\begin{equation}
%    
%    \label{manifoldperts}
%\end{equation}
where $\boldsymbol{\xi_{U,0}}$ and $\boldsymbol{\xi_{S,0}}$ are the unstable/stable eigenvector, $s=\pm 1$ indicates the direction of the manifold, and $\epsilon$ is a small offset. These perturbed initial conditions can be numerically integrated to recover the unstable and stable manifolds in both the forward and backward directions by selecting different values for $s$. Selecting different $\boldsymbol{x}(t_0)$ along the period orbit, recomputing the Monodromy matrix and using and integrating along stable and unstable eigenvalues will yield an approximation of the invariant manifolds associated with a designated periodic orbit (see Fig.~\ref{manifolds}). Invariant manifolds associated with the reference L2 Halo orbit are computed for a short time horizon are shown below. The motivation of this brief discussion on invariant manifolds was to explain the disjoint reachable set and reachable sample trajectory clustering most directly observed in Figure \ref{fig:L2350h}. The unstable and stable manifolds are computed for this same L2 Halo reference periodic orbit and then overlaid onto the reachable set results from Figure \ref{fig:L2350h} to identify a possible trend in reachable trajectories arising from any underlying dynamical structure associated with the chosen periodic orbit. 

 \begin{figure}[!tbp]
\centering 
\subfloat[3D view]{%
  \includegraphics[clip,width=0.8\columnwidth]{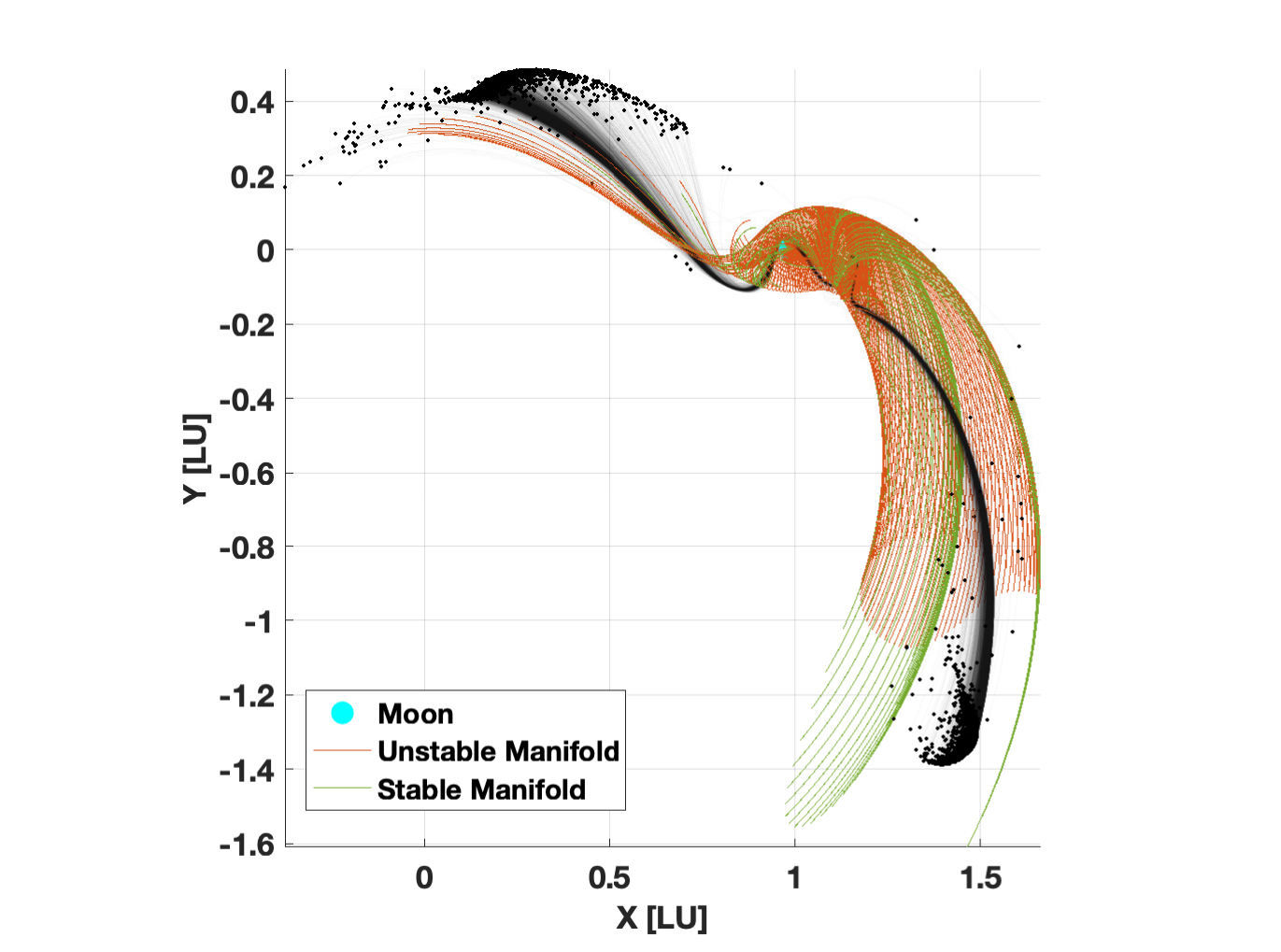}
  \label{fig:reachmanifoldoverlay_xyz}
  }

\subfloat[XY view]{%
  \includegraphics[clip,width=0.8\columnwidth]{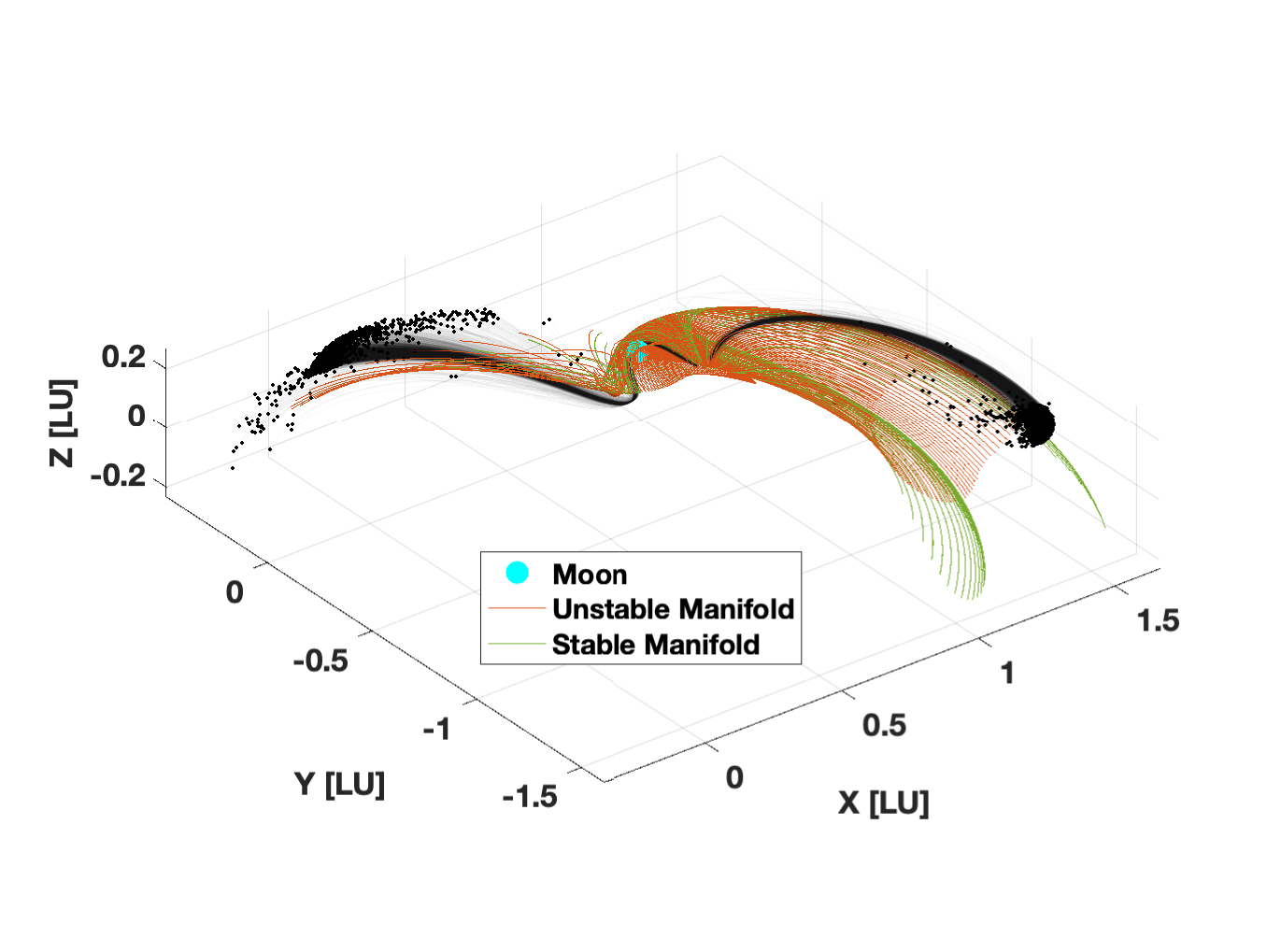}
  \label{fig:eachmanifoldoverlay_xy}
  }
  
\caption{350-hour L2 reference orbit reachable set overlaid with invariant manifolds.}
\label{fig:manifoldoverlay}
\end{figure}

Figure \ref{fig:manifoldoverlay} presents the L2 Halo 350-hour reachable set with the associated invariant manifolds. There is a clear direct relationship, for this long time horizon, of the tendency for reachable set sample trajectories to cluster on the invariant manifolds. Identify that the invariant manifolds propagate in two primary directions originating at the periodic orbit, which is due to the integration in both the forward and backward directions due to variations in the parameter, $s$ \cite{singh_eclipse-conscious_2021}. Since the reachable set is, by definition, the set of all possible states that could be achieved from a given initial condition, it is expected that the sample points terminating the furthest linear distance from the initial condition arrived at the terminal point by exploiting the invariant manifold structure existing in the CR3BP. Additionally, since the propulsion system used for this model (0.2 Newtons) is very low leading to an overall small perturbation acceleration, it is expected that the reachable set evolution on this time horizon is dominated more by the dynamics in the CR3BP than the on-board propulsion system. That is, there is a higher probability of a reachable sample point to terminate near an invariant manifold than any other random spatial coordinate, making these invariant manifolds essentially analogous to supertracks in a general chaotic dynamical system. To our best knowledge, the relation between the reachable set of low-thrust trajectories and invariant manifolds in the CR3BP is not identified in the literature.

\section{Conclusion} \label{sec:conclusion}
The indirect multistage formulation of optimal control theory was used to implement a rapid reachable set approximation algorithm. %The details for modeling a low-thrust spacecraft under two-body and Circular Restricted Three-body Problem (CR3BP) dynamics were presented. 
A minimum-time low-thrust trajectory formulation was used to recover minimum-time trajectories corresponding to reachable sets. 
%This optimal control problem was solved efficiently with the implementation of the reachable set algorithm that uniformly samples terminal costates from a given reference trajectory, backwards integrates to obtain optimal control vectors, and scales and reconstructs a minimum-time reachable set from the sampled trajectories. 
%We detail the rapid computational capabilities of extending time horizons for the Earth-Mars transfer scenario.
Results are presented detailing the application of this algorithm to a minimum-time Earth-Mars transfer. Results indicate that the minimum-time rendezvous lies on the interior of the reachable position set and on the boundary of the reachable velocity set. The effects of injecting initial condition uncertainty and an impulsive maneuver are also evaluated, which shows variations in the reachable set. Novel results detailing the reachable set for low-thrust CR3BP dynamics are presented for L1 point, L2 Halo, and 9:2 NRHO reference conditions were analyzed. Finally, theoretical insight into the behavior of the reachable set, subject to multi-body perturbations, is provided to show that the reachable set and the invariant manifold structures closely align.

The two situations where the proposed algorithm demonstrates areas for improvement are when a reachable set passes close to a gravitational body and in a multi-revolution scenario. For the former, when we observe the reachable set to develop to encompass a gravitational body, the resulting reachable set appears chaotic. The latter is an extension of the long time horizon shortcoming in that for a multi-revolution case the reachable set will overlap with itself and the set becomes incoherent. Our results also indicate that thrust-reversal maneuvers also deform the velocity reachable set. In general, these shortcomings can be avoided completely by computing short time horizon reachable sets, which is the most likely application for a space-domain awareness situation.
%Future work will evaluate the prospect of generating a multidimensional surface from the reachable sample points. This surface will be used for rapid collision avoidance detection and maneuvering. Short time horizon spacecraft proximity operations and satellite constellation operations are also worth investigating in the future.

\section{Appendix: Derivation of Boundary Conditions for An Initial Impulse Maneuver}\label{sec:appendix}

%Here, an  expression for the initial impulse perturbation, $\boldsymbol{\delta v_2}$, due to an impulsive maneuver applied at the initial condition is derived. 
%Recall, the matrix-mapping technique use to augment the initial reference state with perturbed states due to initial condition variations.
%The following derivation will follow a traditional indirect optimal control formulation derivation, that is, define an objective function and constraints, build the Hamiltonian, and form the necessary conditions for optimality to recover expressions for the costates. 
%We will use the same cost functional as the initial position and velocity uncertainty constraints. 
An initial impulsive maneuver introduces the constraint given in Eq.~\eqref{impulsecons}.
%\begin{equation}
%    L^{*}=m_0-m_0 e^{-\frac{\Delta V}{c}}.
%\end{equation}
%\begin{align} \label{impulsecons}
%    \frac{1}{2} \boldsymbol{\delta %v_2}^\top \boldsymbol{\delta v_2}-\frac{1}%{2} \Delta V_{\max }^2 & \leq 0.
%\end{align}
The augmented Hamiltonian is formed as,
\begin{equation}
\begin{aligned}
& H^{*}=m_0-m_0 e^{-\frac{\Delta V}{c}}+\boldsymbol\lambda^\top\left[\boldsymbol{r}^{*^\top},\boldsymbol{v}^{*^\top}\right]^{\top}+ \nu_3\left(  \frac{1}{2} \boldsymbol{\delta v_2}^\top \boldsymbol{\delta v_2}-\frac{1}{2} \Delta V_{\max }^2 \right) . 
\end{aligned}
\end{equation}

We can explicitly define the relationship between $\boldsymbol{\delta v_2}$ and $\Delta V$ as $\Delta V = \left\Vert\boldsymbol{\delta v_2} \right\Vert = \sqrt{\boldsymbol{\delta v_2}^\top \boldsymbol{\delta v_2}}$. Additionally, the derivative of $\Delta V$ with respect to $\boldsymbol{\delta v_2}$ is found as the derivative of a 2-norm function given by,
\begin{equation}
    \frac{\partial \Delta V}{\partial\boldsymbol{\delta  v_2}} = \frac{\partial \left\Vert\boldsymbol{\delta v_2} \right\Vert}{\partial \boldsymbol{ \delta  v_2}} = \frac{\boldsymbol{\delta v_2}}{\left\Vert\boldsymbol{\delta v_2} \right\Vert} =  \frac{\boldsymbol{\delta v_2}}{\Delta V}.
\end{equation}

Forming the necessary conditions for optimality, we obtain,
\begin{align}\label{dHdv2}
     \frac{\partial H^*}{\partial\boldsymbol{\delta v_2} } & = H_{\boldsymbol{\delta v_2}}^*= \textbf{0} =\boldsymbol\lambda_{\boldsymbol{v_2}}+ \frac{\boldsymbol{\delta v_2}}{\Delta V}\frac{m_0}{c}e^{-\frac{\Delta V}{c}} + \nu_3 \boldsymbol{\delta v_2}.
\end{align}

Since Eq.~\ref{impulsecons} is an inequality, the solution process is to assume the constraint is active by setting $\Delta V = \Delta V_{\text{max}}$, obtain an expression for $\nu_3$, and check if $\nu_3>0$. If the inequality is inactive, we will set $v_3=0$ to solve for an expression for $\Delta V$. Applying this logic and substituting into the right relation in Eq.~\ref{dHdv2}, we obtain,
%\begin{equation}
%    \textbf{0} = \boldsymbol\lambda_{\boldsymbol{v_2}}+ \frac{\boldsymbol{\delta v_2}}{\Delta V_{\text{max}}}\frac{m_0}{c}e^{-\frac{\Delta V_{\text{max}}}{c}} + \nu_3 \boldsymbol{\delta v_2} -\lambda_{n}\left(  \frac{\boldsymbol{\delta v_2}}{\Delta V_{\text{max}}} \frac{1}{m_0c} e^{\frac{\Delta V_{\text{max}}}{c}} \right).
%\end{equation}

%Factoring $\boldsymbol{\delta v_2}$, simplifying, and rearranging gives,
\begin{equation} \label{condbindrelate}
    %\textbf{0} = \boldsymbol\lambda_{\boldsymbol{v_2}}+ \left( \frac{m_0}{\Delta V_{\text{max}}c}e^{-\frac{\Delta V_{\text{max}}}{c}} + \nu_3 - \frac{\lambda_{n}}{\Delta V_{\text{max}}m_0c} e^{\frac{\Delta V_{\text{max}}}{c}} \right)\boldsymbol{\delta v_2},
    \boldsymbol\lambda_{\boldsymbol{v_2}} = - \left( \frac{m_0}{\Delta V_{\text{max}}c}e^{-\frac{\Delta V_{\text{max}}}{c}} + \nu_3 \right)\boldsymbol{\delta v_2}.
\end{equation}

%\begin{equation}\label{condbindrelate2}
    
%\end{equation}

In order for the equality condition in Eq.~\ref{condbindrelate} to be satisfied, $\boldsymbol\lambda_{\boldsymbol{v_2}}$ must be parallel to $\boldsymbol{\delta v_2}$. Hence, the optimal solution for $\boldsymbol{\delta v_2}$ when the inequality constraint is binding can be determined from Lawden's Primer Vector theory as,
\begin{equation} \label{primerimpulse}
    \boldsymbol{\delta v_2} = -\Delta V_{\text{max}}\frac{\boldsymbol\lambda_{\boldsymbol{v_2}}}{\left\Vert \boldsymbol\lambda_{\boldsymbol{v_2}} \right\Vert}, \rightarrow \boldsymbol{\delta v_2} + \Delta V_{\text{max}}\frac{\boldsymbol\lambda_{\boldsymbol{v_2}}}{\left\Vert \boldsymbol\lambda_{\boldsymbol{v_2}} \right\Vert} = \textbf{0}, \rightarrow   \boldsymbol\lambda_{\boldsymbol{v_2}} + \frac{\left\Vert \boldsymbol\lambda_{\boldsymbol{v_2}} \right\Vert}{\Delta V_{\text{max}}}\boldsymbol{\delta v_2} = 0.
\end{equation}

%Next, we must obtain an expression for $\nu_3$. Rearranging Eq.~\ref{primerimpulse} to the form of Eq.~\ref{condbindrelate} gives,
%\begin{equation}
   
%\end{equation}

Setting the term in parentheses in Eq.~\eqref{condbindrelate} to $\left\Vert \boldsymbol\lambda_{\boldsymbol{v_2}} \right\Vert / \Delta V_{\text{max}}$,

\begin{equation} \label{nu3_1}
    \frac{\left\Vert \boldsymbol\lambda_{\boldsymbol{v_2}} \right\Vert}{\Delta V_{\text{max}}} =  \frac{m_0}{\Delta V_{\text{max}}c}e^{-\frac{\Delta V_{\text{max}}}{c}} + \nu_3, \rightarrow  \nu_3 = \frac{\left\Vert \boldsymbol\lambda_{\boldsymbol{v_2}} \right\Vert}{\Delta V_{\text{max}}} - \frac{m_0}{\Delta V_{\text{max}}c}e^{-\frac{\Delta V_{\text{max}}}{c}} > 0.
\end{equation}

%Solving for $\nu_3$ and reintroducing the original inequality sign gives,
%\begin{equation} 
%\end{equation}

Only the sign of $\nu_3$ is required, so Eq.~\eqref{nu3_1} can be simplified by multiplying all terms by $\Delta V_{\text{max}}m_0c e^{-\frac{\Delta V_{\text{max}}}{c}}$.
\begin{equation} \label{nu3_2}
    \nu_3 = \left\Vert \boldsymbol\lambda_{\boldsymbol{v_2}} \right\Vert m_0 c e^{-\frac{\Delta V_{\text{max}}}{c}} - \left(m_0 e^{-\frac{\Delta V_{\text{max}}}{c}}\right)^2 > 0.
\end{equation}

Equation \eqref{nu3_2} can be evaluated for each of the sample trajectories in a problem to determine if $\nu_3>0$. If true, Eq.~\eqref{primerimpulse} is the optimal solution for $ \boldsymbol{\delta v_2}$ and $\Delta V = \Delta V_{\text{max}}$. If false, then $\Delta V \neq \Delta V_{\text{max}}$, set $\nu_3=0$ and solve Eq.~\eqref{condbindrelate} for $\Delta V$. Equation \eqref{condbindrelate} with updated notation for a non-binding inequality constraint and $\nu_3= 0$ is 

\begin{equation} \label{condbindnorelate2}
    \textbf{0} = \boldsymbol\lambda_{\boldsymbol{v_2}}- \left( \frac{m_0}{\Delta Vc}e^{-\frac{\Delta V}{c}}  \right)\boldsymbol{\delta v_2}, \rightarrow \textbf{0} = \boldsymbol\lambda_{\boldsymbol{v_2}}- \left( \frac{m_0}{c}e^{-\frac{\Delta V}{c}} \right)\frac{\boldsymbol{\delta v_2}}{\Delta V}.
\end{equation}

We have factored $\Delta V$, since it is related to $\boldsymbol{\delta v_2}$ by $\Delta V = \Vert \boldsymbol{\delta v_2} \Vert$. Similar to the constraint binding case, notice that $ \boldsymbol\lambda_{\boldsymbol{v_2}}$ must be parallel to $\frac{\boldsymbol{\delta v_2}}{\Delta V}$ to satisfy Equation \ref{condbindnorelate2}. That is,
\begin{equation}
    \frac{\boldsymbol{\delta v_2}}{\Delta V} = \frac{\boldsymbol\lambda_{\boldsymbol{v_2}}}{\left\Vert \boldsymbol\lambda_{\boldsymbol{v_2}} \right\Vert}.
\end{equation}

The term in parentheses in Equation \ref{condbindnorelate2} must equal $\left\Vert \boldsymbol\lambda_{\boldsymbol{v_2}} \right\Vert$, which can be written as,
\begin{equation}
    \left\Vert \boldsymbol\lambda_{\boldsymbol{v_2}} \right\Vert = \frac{m_0}{c}e^{-\frac{\Delta V}{c}}, \rightarrow 0 = \frac{m_0}{c}e^{-\frac{\Delta V}{c}} - \left\Vert \boldsymbol\lambda_{\boldsymbol{v_2}} \right\Vert .
\end{equation}

Multiplying by $m_0ce^{-\frac{\Delta V}{c}}$ to simplify and solve the resulting quadratic equation for $m_0e^{-\frac{\Delta V}{c}}$
\begin{equation}\label{quadsol}
    0 = \left( m_0e^{-\frac{\Delta V}{c}} \right)^2 -  \left\Vert \boldsymbol\lambda_{\boldsymbol{v_2}} \right\Vert c \left(m_0e^{-\frac{\Delta V}{c}}\right), \rightarrow m_0e^{-\frac{\Delta V}{c}} = \frac{\left\Vert \boldsymbol\lambda_{\boldsymbol{v_2}} \right\Vert c \pm \sqrt{\left\Vert \boldsymbol\lambda_{\boldsymbol{v_2}} \right\Vert^2 c^2}}{2}.
\end{equation}

To minimize fuel (i.e., maximize the left-hand side of Eq.~\ref{quadsol}), the positive sign has to be selected to have
\begin{equation}\label{deltaV}
    \Delta V =-c \ln \left(\frac{\left\Vert\boldsymbol{\lambda}_{\boldsymbol{v}}\right\Vert c+\sqrt{\left\Vert\boldsymbol{\lambda}_{\boldsymbol{v}}\right\Vert^2 c^2}}{2 m_0}\right), 
\end{equation}

This $\Delta V$ can be used to recover $\delta_n$. The perturbation on initial conditions due to an impulsive maneuver is given by substituting Eq.~\eqref{deltaV} into the primer vector in Eq.~\eqref{primerimpulse},
\begin{equation}
     \boldsymbol{\delta} \boldsymbol{v}_{\boldsymbol{2}} =c \ln \left(\frac{\left\Vert\boldsymbol{\lambda}_{\boldsymbol{v}}\right\Vert c+\sqrt{\left\Vert\boldsymbol{\lambda}_{\boldsymbol{v}}\right\Vert^2 c^2}}{2 m_0}\right) \frac{\lambda_{\boldsymbol{v}}}{\left\Vert\boldsymbol{\lambda}_{\boldsymbol{v}}\right\Vert}.
\end{equation}

\vspace{-5mm}
\section{Disclaimer}
The first author would like to emphasize that the views expressed in this
article are those of the authors and do not reflect the official policy or position of the United States
Air Force, Department of Defense, or the U.S. Government.

\section{Acknowledgments}
The authors would like to thank Dr. Prashant Patel for fruitful discussions about the indirect multistage formulation. 
%We thank Mr. Nick P. Nurre for providing the MEX-friendly version of MATLAB's \texttt{ode113}.
%\vspace{-5mm}
\bibliography{references, misc, taheri}

\end{document}